\input amstex 
\documentstyle{amsppt}  
\input xy
\xyoption{all}
\input graphicx.tex
\magnification=1200
\hsize=150truemm 
\vsize=224.4truemm
\hoffset=4.8truemm
\voffset=12truemm

   \TagsOnRight
\NoBlackBoxes
\NoRunningHeads

\def\Square{\rlap{$\sqcup$}$\sqcap$}
\def\cqfd {\quad \hglue 7pt\par\vskip-\baselineskip\vskip-\parskip
{\rightline{\Square}}}

\define\E{{\Cal E}} 
\define\G{{\Cal G}} 
 
\redefine\S{{\Cal S}} 
\define\T{{\Cal T}} 
\define\M{{\Cal M}} 
\define\C{{\Cal C}} 
\define\A{{\Cal A}} 
\redefine\D{{\Cal D}} 
\define\R{{\Bbb R}} 
\define\Z{{\Bbb Z}} 
\define\Q{{\Bbb Q}} 
\define\Rt {${\R}$-tree} 
\define\mi{^{-1}} 
\let\thm\proclaim 
\let\fthm\endproclaim 
\let\inc\subset  
\let\ds\displaystyle 
\let\ev\emptyset 
\let\wtilde\widetilde
\let\ov\overline
\let\bar\overline
\define\Aut{\text{\rm Aut}}
\define\Out{\text{\rm Out}}
\define\m{^{-1}}
\define\eps{\varepsilon}
\redefine\phi{\varphi}
\define\es{\emptyset}

\define\ra{\rightarrow}
\def\emph#1{{\it #1\/}}
\define\rond#1{\overset{\circ} \to {#1}}

\define\fix{\text{\rm Fix}\,}

\define\normal{\vartriangleleft}

\def\na{non-ascending}

\newcount\tagno
\newcount\secno
\newcount\subsecno
\newcount\stno
\global\subsecno=1
\global\tagno=0
\define\ntag{\global\advance\tagno by 1\tag{\the\tagno}}

\define\sta{\the\secno.\the\stno
\global\advance\stno by 1}

\define\sect{\global\advance\secno by
1\global\subsecno=1\global\stno=1\
\the\secno. }

\define\subsect{\the\secno.\the\subsecno. \global\advance
\subsecno by 1}

\def\nom#1{\edef#1{\the\secno.\the\stno}}
\def\eqnom#1{\edef#1{(\the\tagno)}}

\newcount\refno
\global\refno=0
\def\nextref#1{\global\advance\refno by 1\xdef#1{\the\refno}}
\def\bref {\ref\global\advance\refno by 1\key{\the\refno}}

\nextref\Ba
\nextref\BK
\nextref\BFa
\nextref\Bo
\nextref\BM
\nextref\Chi
\nextref\Cl
\nextref\Clfp
\nextref\Clb
\nextref\Co
\nextref\CM
\nextref\CVa
\nextref\CVb
\nextref\Di
\nextref\DS
\nextref\Fo
\nextref\FoCMH
\nextref\FoGBS
\nextref\FP
\nextref\GL
\nextref\Bar
\nextref\GLt
\nextref\KV
\nextref\LeGD
\nextref\LeGBS

\nextref\McM
\nextref\MS
\nextref\Pa
\nextref\PaENS
\nextref\RS
\nextref\SS
\nextref\Se
\nextref\Ser
\nextref\Sh
\nextref\Sk
\nextref\Wh

\topmatter

  \title  Deformation spaces  of  trees   
            \endtitle

 \abstract     Let $G$ be a finitely generated group. Two
simplicial $G$-trees are said to be in the same deformation space if
they have the same elliptic subgroups (if $H$ fixes a point in one tree,
it also does in the other). Examples include Culler-Vogtmann's outer
space, and spaces of JSJ decompositions. We discuss what features
are common to trees in a given deformation space, how to pass from
one tree to all other trees in its deformation space, and the topology
of deformation spaces. In particular, we prove that all deformation
spaces are contractible complexes.
\endabstract

\author  Vincent Guirardel, Gilbert Levitt
        \endauthor

\endtopmatter

\document

\head     Introduction\endhead 

Let $G$ be a finitely generated group. A $G$-tree is a simplicial tree
with an action of $G$. The notion of a deformation space was introduced
by Forester [\Fo]. By  definition, two $G$-trees are in the same
deformation space
$\D$ if they have the same elliptic subgroups: if $H\inc G$ fixes a point
in one tree, it also does in the other. Forester showed that two trees in
the same $\D$ may be connected by a finite sequence of elementary
deformations, each associated to a  canonical isomorphism $A*_B
B\simeq A$. 

Examples of deformation spaces are Culler-Vogtmann's outer space
 [\CVa], as well as spaces constructed by McCullough-Miller [\McM]
and the authors [\GL] to study automorphisms of free products, and
the canonical set of   splittings of generalized Baumslag-Solitar
groups [\FoGBS,
\LeGBS].

Deformation spaces are especially relevant to JSJ theory. JSJ
decompositions of finitely presented groups have been constructed by
Rips-Sela [\RS], Dunwoody-Sageev [\DS], Fujiwara-Papasoglu [\FP].
They are
$G$-trees (equivalently, graphs  of groups decompositions of $G$)
with certain properties. Though canonical, they are only unique up to
certain moves (see the above references, as well as [\FoCMH], [\Bo],
[\SS]). 

We shall explain in [\GLt]   (see [\Bar]) that, in general, the
canonical object is not a JSJ-tree but a JSJ deformation space, and
that a JSJ-tree is just as unique as a tree in a general  deformation
space.  We shall also give a general construction of this
JSJ  space, valid for any finitely presented group and any class of edge
groups (not just slender ones).

With this in mind, we study here  general properties of
deformation spaces, focusing on three main questions:

$\bullet$ Starting with one tree, what moves are needed to generate all
trees in its deformation space? Are slide   moves sufficient?

$\bullet$  What is common to trees in the same deformation space? In
particular, to what extent do trees in the same space have the same
vertex and edge stabilizers?

$\bullet$  What is the topology of a deformation space? Is it
contractible? Is it finite-dimensional, or with a finite-dimensional
spine?

\head   \sect Contents of the paper\endhead

We fix a finitely generated group $G$ and we consider metric simplicial
trees $T$ with an isometric action of $G$, up to equivariant isometry.
The action is always assumed to be minimal. In this introduction we
also assume irreducibility (it makes several statements simpler).

A subgroup $H\inc G$ is elliptic in $T$ if it fixes a point. By
definition, two trees
$T,T'$ are in the same deformation space $\D$ if they have the same
elliptic subgroups. Equivalent characterizations (mainly due to
Forester [\Fo]) are the following (see Theorem   3.8): there exist
equivariant maps $f:T\to T'$ and $f':T'\to T$; the trees are related by a
sequence of elementary deformations; their length functions are
bi-Lipschitz equivalent. 

We also consider
restricted deformation spaces $\D_\A$, where we only consider trees
with edge stabilizers in a class $\A$ of subgroups of $G$ (such as
cyclic subgroups, abelian subgroups, etc.).

Sections 2  and 3  consist mainly of known facts and definitions. In
Section 4, we discuss what invariants may be extracted from a given
deformation space $\D$.

 For instance, the first Betti
number of the quotient graph $\Gamma =T/G$ depends only on $\D$, not
on $T$. We  denote  it by $b_1(\D)$. 

Unfortunately, it is not true that all trees in  
$\D$ have the same vertex or edge stabilizers. 
For  one thing, one may
(almost always) introduce new vertex groups by using the
isomorphism 
$A\simeq A*_B B$. This  leads us to consider  trees which are reduced
(in the sense of [\Fo]): if one collapses to a point all edges in a given
orbit, one always obtains a tree  outside of $\D$. 

But reduced trees in
$\D$ don't always have the same vertex groups. A basic example is
the Baumslag-Solitar group $BS(1,6)=\langle a,t\mid
tat\mi=a^6\rangle$ [\FoGBS]. In the deformation space  
consisting of trees with all stabilizers
cyclic, there exists a reduced
  tree with vertex stabilizers conjugate to $
\langle a \rangle$, and another one with stabilizers conjugate to
$\langle a^2\rangle$. If $G=BS(2,4)$, there exists a reduced   tree
with one orbit of vertices, and another tree with two orbits.

One way to avoid these problems is to consider a restricted
deformation space   $\D_\A$ and vertex stabilizers not in $\A$. It is then
true that all trees in  $\D_\A$ have the same vertex stabilizers not in
$\A$. 

 In general, the information about vertex
stabilizers of trees in
$\D$ is captured  by a finite set $\M$ which we associate to $\D$. 

In nice
situations, for instance when no group fixing an edge is properly
  contained in a conjugate of itself,
$\M$ is simply the set of conjugacy classes of vertex stabilizers of
reduced trees (see
Section 7 for more general statements). In general, an element of
$\M$ is a union of conjugacy classes of elliptic subgroups.

 If $\T\in\D$ is
reduced, we   show that the  number $s$  of vertices of the quotient
graph  satisfies
$|\M|\le s\le 2|\M|+2b_1(\D)-2$. 
This   upper bound on $s$ can be viewed as a simple accessibility
result, holding within a fixed deformation space (note that $G$ is not
assumed to be finitely presented).

Similarly, it is not always true that reduced trees in $\D$ have the
same edge groups. But they have the same  
bi-elliptic groups, and the same generalized edge groups, where
a group is bi-elliptic if it fixes two distinct points of $T$, is a
generalized edge group if it is bi-elliptic and furthermore it contains
some edge stabilizer.

In JSJ theory, one encounters quadratically hanging subgroups. The
important information about such a subgroup is not its isomorphism
type (it is often a free group), but the way it is attached to the rest of
the group (the topological picture being a compact surface attached
along its boundary components). It is therefore desirable to attach a
peripheral  structure to vertex stabilizers of  trees in $\D$. 

We explain how to define
such a structure, given  $\D$  and a maximal elliptic subgroup 
$G_0$. All trees in
$\D$ have a vertex $v$ with stabilizer $G_0$, and the peripheral
structure   contains information about incident edge groups. In nice
situations, it is simply the set of  conjugacy classes of stabilizers of
edges incident to
$v$.

In Section 5, we discuss topologies on a deformation space $\D$. There
are two natural topologies on $\D$. The first topology is the
equivariant Gromov-Hausdorff topology (or simply Gromov topology),
also called the axes topology because it may be described in terms of
length functions. The second topology is the weak topology,
associated to the natural structure of a cell complex on $\D$. 

As 
pointed out in   [\McM], these topologies may be different. We show, 
however,  that they always agree on any  subset of $\D$ which is
contained in a finite union of cells, and also that they agree on the
whole of $\D$ when $\D$ consists of locally finite trees with finitely
generated stabilizers. (As mentioned above, we assume in this
introduction that trees are irreducible. When
$\D$ consists of $G$-trees with a  $G$-fixed end, we have to show
the non-obvious fact that $\D$ is Hausdorff in the Gromov topology). 

In Section 6, we   discuss contractibility.  Skora [\Sk]  introduced the 
idea of deforming morphisms to prove that Culler-Vogtmann space  
and its closure are contractible. Using his technique, we show that a
deformation space
$\D$ is always contractible in the weak topology. When $\D$ contains a
tree with finitely generated vertex groups, we show that it is
contractible in the Gromov topology (this was proved independently
by Clay [\Cl]), and also that the closure of $\D$ is contractible. 

  Here is a sketch of the argument for proving contractibility. It is
relatively easy to construct a map
$F:\D\times[0,1]\to\D$
 contracting $\D$ to a point (this requires 
choosing a
  basepoint in trees of $\D$  in a  continuous way; we give a direct
geometric construction, different from Skora's minimization
argument).
 The problem is to show that $F$ is continuous. We
consider its restriction   to
$\S\times[0,1]$, where
$\S$ is a   closed cell. It is  continuous in the Gromov topology. We
show that the image of this restriction  is  a subset of $\D$ meeting
only  finitely many cells (finiteness lemma 6.4). Since  the two
topologies agree on such a subset, $F$  is continuous in the weak
topology on $\S\times[0,1]$, hence on the whole of $\D\times[0,1]$
by definition of the weak topology. As in [\Cl], proving continuity in the
Gromov topology requires a finiteness hypothesis.

We also use   Skora's  technique to give a direct proof of Forester's
deformation theorem [\Fo].

In Section 7, we study a ``nice'' class of deformation spaces, the \na{}
ones. Without giving the definition here, let us mention two special
cases. First, spaces containing a tree which is weakly acylindrical
in the following sense: if $g\in G$ is nontrivial, its fixed point set
contains no infinite ray. Also, spaces consisting of trees such that no
edge stabilizer properly contains a conjugate of itself. This applies in
particular to splittings over finite subgroups, and to cyclic splittings
of groups containing no solvable Baumslag-Solitar group $BS(1,n)$ with
$n\ge2$. 

Generalizing a result of Forester, we show that any two reduced trees
in a \na{} deformation space are related by slide moves. In particular,
they have  exactly the same vertex and edge stabilizers. We also show
that a \na{} deformation space has a natural
deformation retraction onto a finite-dimensional subcomplex.

In Section 8, we consider a finitely generated subgroup $F$ of $\Out(G)$
leaving a deformation space $\D$ invariant (see  [\BM, \Clfp, \KV, \Wh]
for related results). We show that the fixed point set of $F$ in $\D$ is
empty or contractible, using the fact that this fixed point set is a
deformation space  of
$\widehat G$-trees for some extension $\widehat G$ of $G$. If $F$ is a
solvable finite group and 
$\D$ is \na{}, we show that $F$ does have a fixed point in $\D$. We
conclude the paper by a few facts   about automorphisms leaving
invariant a deformation space consisting of locally finite trees.

\head   \sect Trees\endhead

Let $G$ be a finitely generated group. Unless otherwise indicated (in
 Sections 5 and 6),  all trees   will be simplicial
$G$-trees, i.e\. we consider   a simplicial tree $T$ with an action of $G$
by simplicial automorphisms, without inversions. See  [\Chi,
\CM, \Sh]  for   basic facts about trees.

We always assume
that
$T$ is {\it minimal\/}: there is no proper $G$-invariant subtree. This
implies that there are only finitely many orbits of edges, and no
terminal  vertices. 
 
We usually   assume
that there is no redundant vertex (every vertex of $T$ has degree
$\ge3$), though we sometimes need to subdivide trees (in the proof of
Lemma 6.5, for instance).

Two distinct points $a,b$ of $T$ bound a unique  {\it segment\/} $[a,b]$.
A {\it finite subtree\/}  is the convex hull of a finite set of points.  An
{\it end\/} of $T$ is an equivalence class of infinite rays, two rays being
equivalent if their intersection is a subray.

A       $G$-tree $T$ may be considered as a simplicial  (or
combinatorial) object, or as a {\it metric space\/} with an isometric
action of $G$. The simplest way to define a metric is to assign length 1
to every edge. More generally, one may assign any positive length to
each orbit of edges (recall that there are finitely many orbits). 

Different length assignments lead to metric trees with the same
underlying simplicial tree (in Section 5, we will say that they belong to
the same open cone). They are equivariantly bi-Lipschitz
homeomorphic. It will sometimes be important to distinguish between
a metric tree and the underlying   simplicial (non-metric) tree.

Let $T$ be a   metric $G$-tree.
The {\it length function\/}   $\ell_T:G\to\R$ is defined by
$\ell_T(g)=\min_{x\in T}d(x,gx)$. When no confusion is possible, we 
simply write $\ell$. An element
$g\in G$ is {\it elliptic\/} if it has a fixed point (equivalently, if  
$\ell(g)=0$), {\it hyperbolic\/} otherwise. A hyperbolic element has a
{\it translation axis\/} $A_g$, on which it acts as translation by
$\ell(g)$, and    
$d(x,gx)=\ell(g)+2d(x,A_g)$ for   any $x\in T$. This formula also
holds if $g$ is elliptic, with $A_g$ understood as the fixed point set of
$g$.

A subgroup $H\inc G$ is {\it elliptic\/} if it fixes a point. A finitely
generated subgroup is elliptic if and only if all its elements are
elliptic [\Ser]. But an infinitely generated subgroup consisting of 
elliptic elements may fail to be elliptic; it then fixes a unique end of
$T$  (see for instance [\Chi, Theorem 3.2.6]). 

We distinguish five {\it types\/} of minimal $G$-trees (see   [\CM] or
[\Chi{}, p\. 134]):

$\bullet$ {\it Trivial\/}: $T$ is a point.

$\bullet$ {\it Dihedral\/}: $T$ is a line, but $G$ does not preserve
orientation. The action factors through an action of the infinite
dihedral group $D_\infty\simeq\Z/2      \Z*     \Z/2    \Z$.

$\bullet$ {\it Linear abelian\/}: $T$ is a line, and $G$ acts by translations.
The action  factors through an action of $\Z$. 

$\bullet$ {\it Genuine abelian\/}: $G$ fixes an end of $T$, and $T$ is not a
line. The quotient  graph $\Gamma =T/G$ is homeomorphic   to a
circle. When there is only one orbit of edges, $T$ is  the Bass-Serre
tree of a strictly ascending HNN-extension $G=\langle A,t\mid
tat\mi=\varphi (a)\rangle$, where
$\varphi :A\to A$ is injective but not onto.  

$\bullet$ {\it Irreducible\/}: there exist two hyperbolic elements with
disjoint axes. In this case, $G$ contains a free group of rank 2 acting
freely.    

Knowing the partition of $G$ into elliptic and hyperbolic
elements (i.e\. the zero-set of the length function) is
enough to distinguish   trivial, dihedral,
abelian, irreducible trees [\CM] (but the two abelian types are not distinguished). In
particular,
$T$ is irreducible if and only if there exist hyperbolic elements
$g,h$   whose commutator    $[g,h]=ghg\mi h\mi$ is hyperbolic. 

A tree is abelian if and only if its length function is the absolute
value of a homomorphism from $G$ to $\R$ with cyclic image. Every element of
the commutator subgroup $[G,G]$ is elliptic, but the subgroup
$[G,G]$ is elliptic if and only if the tree is linear abelian. In
particular, there is no genuine abelian $G$-tree if $[G,G]$ is finitely
generated.

Trees are considered equal if they are equivariantly isometric.  
Two minimal non-abelian $G$-trees with the same length function  are
equivariantly isometric (see  [\CM]).

Maps between $G$-trees will always be $G$-equivariant. As trees
are minimal, maps are always onto. A
map $f:T\to T'$  between simplicial trees  is {\it simplicial\/} if each
edge is mapped bijectively   onto an edge  (in particular, no edge is
collapsed).  It  is a {\it morphism\/}  if
each edge of $T$ can be written as a  finite union of
subsegments, each of which is mapped
bijectively onto a segment in $T'$. Equivalently, $f$ is a
morphism if and only if one may subdivide
$T$ and $T'$ so that $f$  becomes simplicial. When $T$ and $T'$ are
metric trees, simplicial maps and morphisms are required to map  the
subsegments isometrically.

\head   \sect Generalities on deformation spaces  \endhead

We fix a finitely generated group $G$. All trees considered in this 
section will be minimal metric simplicial $G$-trees, up to equivariant
isometry. We always assume that $G$ acts without inversions.  
In this section, and in the next
one, only the simplicial structure  will really   matter.
In later sections, however, allowing edge lengths to vary will
be crucial.

\definition{Definition \sta{} (domination)}
  We say   that $T$ \emph{dominates} $T'$ if there is an equivariant (but
otherwise arbitrary) map
$f$ from $T$ to $T'$. 
\enddefinition

Equivalently, $T$ dominates $T'$ if and only if 
every elliptic subgroup of $T$ is elliptic in $T'$.
The ``only if'' direction   is clear and 
the converse is proved by  first defining $f$ equivariantly on vertices,
and then extending   linearly to edges.
Domination is a pre-ordering. The corresponding equivalence classes 
are \emph{deformation spaces}:

\definition{Definition \sta{} (deformation space)}
The    {\it deformation space\/} $\D$ containing $T$ is the
set of  metric  trees $T'$ such that $T$ and $T'$ dominate each other,
up to equivariant isometry. Trees are in the same $\D$ if and
only if they have the same elliptic   subgroups. Identifying two
trees when they differ only by rescaling the metric leads to the {\it
projectivized deformation space
$P\D$\/}.
\enddefinition

For now $\D$ and $P\D$ are just   sets, with no extra structure, but in 
Section 5 we will view them as topological spaces and complexes.

\definition{Definition \sta{} (collapse, collapsible   edge)}
Let $e$ be an edge of $T$.
 Collapsing every edge in the orbit of
$e$ to  a point produces a new tree $T'$. We  say that $T'$ is
obtained by {\it collapsing\/}    $T$ (or $e$), and $T$ is obtained   by
{\it expanding\/}  $T'$. Collapsing  is defined for simplicial trees. 
In the context of metric
trees,  it   may be viewed as assigning length
$0$ to edges in the orbit of $e$.

The     tree $T'$ is dominated by $T$. Even if $T'$ is not   trivial, it is 
not necessarily in the same
$\D$ as $T$ (this may happen   if $T$ is   irreducible).  We say that $e$
is
\emph{collapsible} if
$T$ and
$T'$ are in the same deformation space.
\enddefinition

The stabilizer of a point $v$, or an edge $e=vw$, will be denoted by 
$G_v$ or
$G_e$. 
We note that an edge $e=vw$ is collapsible if and only if $v$ and $w$
are in   distinct orbits and at least one of the inclusions from $G_e$ to
$G_v$ or from
$G_e$ to
$G_w$ is onto [\Fo]. 

Indeed, if $v$ and $w$ are in distinct orbits and   the inclusion from
$G_e$ to $G_v$ is onto, then collapsing $e$ does not change the set of
elliptic subgroups: the stabilizer of $w$ does not change (since
$G_v*_{G_v}G_w=G_w$), and although  the vertex 
$v$ ``disappears'' its stabilizer remains elliptic (it is contained in
$G_w$).  On the other hand,   collapsing $e$ with $G_e\subsetneqq G_v$
and $G_e\subsetneqq G_w$ would create  a new elliptic subgroup,
namely $G_v*_{G_e} G_w$.  If $v$ and $w$ are in the same orbit, say
$gv=w$, then $g$ is hyperbolic in $T$  (because there is no inversion)
and collapsing $e$ makes $g$ elliptic.

\definition{Definition \sta{} (elementary   collapse, elementary
deformation)} When $e$ is collapsible, passing from $T$ to $T'$ is
called an {\it elementary collapse\/}.    
The reverse  of an elementary
collapse is an {\it  elementary  expansion\/}. A finite sequence of
elementary collapses and expansions  is called an {\it
elementary deformation\/}  [\Fo]. An elementary
deformation between metric trees is an elementary deformation 
between the  underlying simplicial trees. 
\enddefinition

\nom\BFr
\definition{Definition \sta{} (reduced)} 
A tree $T$ is \emph{reduced} (in the sense of Forester
[\Fo]) if no elementary  collapse is possible (collapsing any   edge
yields a tree in a different deformation space).
Equivalently: if an edge $e=vw$ satisfies $G_e=G_v$, then $v$ and $w$
are in the same $G$-orbit.
\enddefinition

Being reduced  is stronger than being
{\it BF-reduced\/} (in the sense of Bestvina-Feighn [\BFa]), where
$G_e\to G_v$ being onto only implies that the image of $v$ in the
quotient graph
$\Gamma =T/G$ has valence at least 3.

Given any tree $T$, one may  
perform a (usually non-unique) sequence of
elementary   collapses so as to obtain a reduced tree $T_r$. We say
that $T_r$ is a {\it reduction\/} of $T$.  

\nom\passage
\example{Remark \sta } 
If the  reduced trees $T,T'$ are related by an elementary
deformation, there exists a sequence of reduced trees $T=T_1,\dots,
T_k=T'$ such that
$T_i$ and $T_{i+1}$ are different reductions of a tree $S_i$ (to see
this, call $S_i$ the   intermediate trees appearing in
the deformation between $T$ and $T'$, and choose a reduction $T_i$ for
each $S_i$). 
\endexample

\nom\pasecr
\example{Remark \sta} Suppose $T,T'$ have the same elliptic
subgroups, and $T$ is reduced. If $f:T\to T'$ is equivariant, and $e=vw$
is an edge of $T$, then
$f(v)\neq f(w)$. 
  Otherwise, one could modify $f$ to be constant on $e$, so $f$ would
factor through the tree $T''$ obtained by  collapsing $e$. But since $T$
is reduced,
$T''$ does not dominate $T$. Since $T''$ dominates $T'$ which dominates
$T$, this is a contradiction.
\endexample

By definition, trees related by an elementary deformation are in the 
same deformation space.  Forester proved that the converse is true,
hence   the terminology.  More explicitly:

\nom\lessix 
\thm {Theorem \sta} Given two $G$-trees $T$ and $T'$, the
following are equivalent:
\roster
\item $T$ and $T'$ are in the same deformation space;
\item there exist equivariant maps $f:T\to T'$ and $f':T'\to T$; 
\item $T$ and $T'$ have the same elliptic subgroups;
\item $T$ and $T'$ are related by an elementary deformation;
\item there exists an equivariant quasi-isometry from $T$ to $T'$;
\item there exists an equivariant map $f:T\to T'$ having bounded
preimages.
\endroster
If $T$ and $T'$ are not abelian, one may add:
\roster
\item[7] the length functions $\ell_T$, $\ell_{T'}$ of $T$ and $T'$ are
bi-Lipschitz equivalent (there exists $C>0$ such that
$\frac1C\ell_T\le\ell_{T'}\le C\ell_T$).
\endroster
\fthm

  All
conditions are obviously satisfied when $T,T'$ are 
metric trees with the same underlying simplicial tree, so this is really 
a statement about simplicial trees.

\demo{Proof} The equivalence between (1) through (6) is in [\Fo].
We sketch some of the arguments.
  We have seen  $(1)\Leftrightarrow(2) 
\Leftrightarrow  (3)$.  
The implication  $(3) \Rightarrow  (4)$ is one of the main
results of [\Fo]; we shall provide a new proof in Section 6 (see 
Theorem 6.6). $(4) \Rightarrow  (5)$ because an elementary
 collapse is a quasi-isometry: if
$e$ is a collapsible edge, then any connected component  of $G.e$ has
diameter at most 2. Obviously $(5) \Rightarrow  (6)$. We show  $(6)
\Rightarrow  (3)$. If $H\inc G$ fixes   $x\in T$,  it fixes $f(x) \in T'$.
Conversely, if
$H$ fixes
$y\in T'$, it leaves invariant   the bounded set $f\mi(y)\inc T$ so   fixes
a point in $T$. Thus $(6) \Rightarrow  (3)$.

We now show $(5)\Leftrightarrow(7)$.  
Given any  $x\in T$ and $g\in G$, one
has
$\ds\ell_T(g)=\lim_{n\to\infty}\frac{d(x,g^nx)}n$, so $(5)$ always
implies $(7)$. Conversely, assume $(7)$. Then $T$ and $T'$ have the
same elliptic elements. We must show that they have the same
elliptic subgroups, provided that they are not abelian.  We suppose 
that a subgroup $H$ fixes
$y\in T'$ but is not elliptic in $T$, and we obtain a contradiction. 

The group    $H$ fixes a unique end of $T$.  Since
$T$ is not abelian, it is irreducible and we can find a hyperbolic
$g\in G$ whose axis $A_g$ does not contain that end. The set
$\ell_T(gH)$ is then unbounded, as there exist elements of $H$ whose
fixed point set is arbitrarily far from $A_g$ (see [\Pa, Proposition
1.7]). By $(7)$, the set
$\ell_{T'}(gH)$ also is unbounded. This is a contradiction, as it is
bounded by the distance from $y$ to
$gy$ in $T'$.  
\cqfd\enddemo

\nom\bdd
\example{Remark \sta} $(3) \Rightarrow(5)$  and  $(3)\Rightarrow(6)$  
may obviously be strengthened as follows:  if
$T,T'$ have the same elliptic subgroups,   then any equivariant map
$f:T\to T'$  which is piecewise linear on each edge
is a quasi-isometry and has bounded preimages. We give a
direct argument to show  that   $f$ has bounded preimages if (3)  
holds (this will be useful in the proof of Theorem 6.6). 
One can subdivide $T$ and $T'$ so that $f$ maps each edge to a vertex 
or an edge. We only need to prove that the preimage of a vertex is
bounded, and  since
$T$, $T'$ have finitely many orbits of   vertices, it suffices to show
that sets of the form
$G.v\cap f\mi(f(v))$ are bounded.  This intersection is equal to
$H.v\cap f\mi(f(v))$, where
$H$ is the stabilizer of $f(v)$. By $(3)$, $H$ is elliptic in $T$, so  $H.v$ is
bounded.
\endexample

\nom\poin
\thm{Proposition \sta}\roster
\item All trees in a given deformation space have the same type
(trivial, dihedral, trivial abelian, genuine abelian, irreducible).
\item All trees in a given abelian deformation space have the same 
length function up to scaling.  
\item If $T$ is trivial,   linear abelian, or dihedral, its
projectivized deformation space consists of only one point. 
\endroster
\fthm

A deformation space $\D$ will be called {\it genuine abelian\/}, or
{\it irreducible\/}, if the trees contained in $\D$ are genuine abelian,
or irreducible. These are the only interesting types.

\example{Remark \sta}  
More generally, a non-elementary collapse does not change the type,
as long as the collapsed tree is not trivial. 
\endexample

\demo{Proof}  The first assertion follows
from remarks  made in Section 2.   

The  length function of an
  abelian   tree is   the
absolute value of a homomorphism   $\varphi :G\to  \R$ with cyclic
image. The deformation space determines the kernel of $\varphi $,
hence
$\varphi $ up to scaling. 
This shows (2).

(3) is trivial in the trivial case.  There is only one  linear abelian   tree
with   length function $|\varphi |$, so (3) holds in the linear abelian case.
Now suppose that  $T$ is dihedral, and $T'$ is in the same deformation
space. The actions of $G$ on $T$ and $T'$ have the same kernel
(elements acting
  as the identity), equal to the set  of   elliptic elements whose product
with any elliptic element is elliptic. One then observes that the infinite
dihedral group only has one non-trivial action, up to scaling. 
\cqfd\enddemo

A deformation space  $\D$  is determined by the family $\E$
consisting of all elliptic subgroups.    Often one also wants to
restrict edge stabilizers. Let $\A$ be a family  of subgroups of $G$,
stable under conjugating and passing to a subgroup. For instance,
$\A$ may consist of all subgroups which are finite,   cyclic,  
abelian, slender, small...

\nom\restri
\definition{Definition \sta{} (restricted deformation space $\D_\A$)}
We define $\D_\A\inc \D$ as the set of trees $T\in \D$ whose edge
stabilizers belong  to $\A$. We call $\D_\A\inc \D$ a {\it restricted
deformation space\/}, or simply a deformation space. We shall see
later that any $\D$ contains a smallest nonempty $\D_\A$, called
the \emph{reduced deformation space} $\D_r$.
\enddefinition

\example{Examples} $\bullet$ If $G$ is a free group $F_n$,
Culler-Vogtmann's outer space is the projectivized deformation
space obtained when $\E$ consists only of the trivial group. 

$\bullet$ Suppose
$G=G_1*\dots*G_p*\Z*\dots*\Z$, where each $G_i$ is non-trivial,
non-isomorphic to $\Z$, and freely indecomposable. Let $\E$ consist
of subgroups of conjugates of the factors, and $\A$ consist only of
the trivial group. The projectivized space $P\D_\A$ is the space
constructed by McCullough-Miller [\McM]. If $\E$ contains only
subgroups of conjugates of the factors $G_i$, then $P\D_\A$ is the 
space considered in [\GL].

$\bullet$ Let $G$ be a generalized Baumslag-Solitar group (GBS group),
i.e\. $G$ acts on a tree with all edge and vertex stabilizers isomorphic
to $\Z$. If
$G$ is not $\Z$, $\Z^2$, or the Klein bottle group,   all such $G$-trees
belong to the same  deformation space [\Fo]. We shall call it the {\it
canonical deformation space\/} of $G$. It is genuine abelian if
$G$ is  a solvable Baumslag-Solitar group $BS(1,n)$ with $n\neq\pm1$,
irreducible if not. If $G=\Z^2$, there are infinitely many deformation
spaces, all of them linear abelian. If $G$ is the Klein bottle group, there
are two spaces; one is linear abelian, the other one is dihedral.
 
$\bullet$ Let $G$ be the fundamental group of a closed orientable
surface $\Sigma
$. Any $G$-tree with cyclic edge stabilizers is dual to a family of
disjoint essential simple closed curves (see   [\MS, Theorem III.2.6]).
After projectivization, the space of metric
$G$-trees with  cyclic edge stabilizers may be identified with the
curve complex of $\Sigma
$. Two such trees are in the same  deformation space if and only if they
belong to the same open simplex. 

$\bullet$ Suppose $G$ has a nontrivial center $Z$, and $\D$ is genuine
abelian or irreducible. Standard arguments show that $Z $ acts as the
identity on any $T\in\D$ (an element $z\in Z$ cannot be hyperbolic, as
its axis would be an invariant line; so
$z$ is elliptic, and its fixed point set is an invariant subtree, equal to
$T$ by minimality). We may therefore view $\D$ as a deformation space
over $G/Z$. 
 
\endexample

\head   \sect Invariants of deformation spaces\endhead

Trees in a given deformation space $\D$ may be fairly different.
Still, certain features depend only on $\D$.  

\subhead Betti number \endsubhead

The quotient of $G$ by the subgroup generated by all elliptic
elements is a free group, isomorphic to the (topological)
fundamental group of the quotient graph $\Gamma =T/G$ for any
$T\in \D$. The rank of this group will be called the {\it Betti
number\/} $b_1(\D)$ of
$\D$.

\subhead Acylindricity\endsubhead

A $G$-tree $T$ is called {\it acylindrical\/} [\Se] if there exists $k$ such
that, for any nontrivial $g\in G$, the fixed point set of $g$ has
diameter less than $k$. As all trees in $\D$ are equivariantly
quasi-isometric to each other, acylindricity depends only on $\D$
(but $k$ depends on $T$). 

Similar considerations apply to the  following weak form of
acylindricity: if $g$ is nontrivial, its fixed point set contains no
infinite ray. 

\subhead Local finiteness\endsubhead

It is easy to see that an elementary collapse, or an elementary
expansion, preserves the property that $T$ is locally finite (recall
that we only consider minimal  trees).   Local finiteness of trees
is therefore a property of the deformation space. For instance
Culler-Vogtmann space, the space constructed by McCullough-Miller
when
$G$ is a free product of finite groups, the canonical deformation space
of a GBS group, have this property.

\subhead Vertex groups\endsubhead

Though all trees in a deformation space $\D$ have the same elliptic
subgroups, they do not necessarily have the same vertex
and edge stabilizers, even if they are reduced. 

\nom\bizz
\example{Examples \sta} 
$\bullet$ In the canonical deformation space of 
$BS(1,6)=\langle a,t\mid tat\mi=a^6\rangle$, there exists a tree $T_1$
whose vertex stabilizers   are conjugates of $\langle a \rangle $, and
another tree $T_2$ whose vertex stabilizers   are conjugates of
$\langle a^2
\rangle $ (see Example 4.3 in [\FoGBS]). Note that $\langle a 
\rangle $ is conjugate to the  subgroup $\langle a^6
\rangle \inc\langle a^2
\rangle $.

$\bullet$ Now consider the Bass-Serre tree of the HNN-extension
$BS(2,4) =\langle a,t\mid ta^2t\mi=a^4\rangle$.   There is one
orbit of vertices. Vertex stabilizers are conjugates of 
$\langle a 
\rangle $. On the other hand, the tree associated to the
presentation
$BS(2,4) =\langle a,b,t\mid tbt\mi=b^2, a^2=b^2\rangle$, though
reduced,  has two orbits of vertices, with  vertex
stabilizers conjugates of 
$\langle b\rangle$ and  $\langle a\rangle$  (see [\LeGBS]).  \endexample

These examples are possible because there exist edge stabilizers 
which properly contain conjugates of themselves.  In Section 7, we will
study \na{} deformation spaces, where such examples cannot occur. In a
\na{}  space $\D$, every elliptic subgroup is contained in a maximal one.
All reduced trees   $T\in\D$ have the same vertex and edge stabilizers,
and there is a natural   bijection  between the set of  $G$-orbits of
 vertices    and the set $\M$ of conjugacy classes of maximal
elliptic subgroups (but there may exist distinct vertices of $T$ with 
the same stabilizer). 

In a general deformation
space, we have to define a set $\M$ in a more complicated way.

\definition{Definition \sta{} (vertical set $\M$)} We associate  a set
$\M$ to $\D$ as follows (compare  the definition of vertical subgroups  
in [\FoGBS]). Given   elliptic subgroups $H$ and $K$, define $H\le K$ if $H$
is contained in a conjugate of $K$, and $H\sim K$ if $H\le K$ and $K\le H$.
The set $\C$ of equivalence classes is partially ordered, and we define
$\M$ as the set of maximal elements of $\C$. 
\enddefinition

If $H$ is a maximal elliptic subgroup, its
equivalence class equals its conjugacy class and is an element
of
$\M$; in general, an element of
$\M$ is a union   of conjugacy classes of elliptic
subgroups. 
If $\D$ contains  $T$   such that no edge-to-vertex
inclusion
$G_e\to G_v$ is onto,   $\M$ is simply the set of conjugacy
classes of maximal elliptic subgroups; it  may   be identified with  the  
set
 of $G$-orbits  of vertices. 
 
\example{Example} 
In the example of $BS(1,6)$ given above, the groups $\langle
a 
\rangle $ and  $\langle
a^2
\rangle $ are equivalent. Their class is the unique element of $\M$.
There is no maximal elliptic subgroup (in this example  $\D$ is
abelian, but this  may also occur in irreducible
deformation spaces). 
In the
$BS(2,4)  $ example, the set $\M$ has one element, the
conjugacy class of $\langle a 
\rangle $ (a maximal elliptic subgroup).  The vertex
stabilizer
$\langle b\rangle$ is not maximal elliptic. 
\endexample

Given a tree $T$, we  denote by $s$ the number of vertices of the
quotient graph $\Gamma =T/G$. It is the number of $G$-orbits of
vertices of $T$. 

\nom\findim
\thm{Proposition \sta} Let $\D$ be a deformation space.
\roster
\item The set $\M$ is finite. If $T\in\D$, the number $s$ of 
 vertices of $\Gamma =T/G$ satisfies
$s\ge |\M|$. 
\item If $T\in\D$ is reduced, or only BF-reduced, then $  s\le
2|\M|+2b_1(\D)-2$. 
\endroster
\fthm

\demo{Proof} Fix $T\in\D$.  There is a natural map  $\varphi
$ from the
set of
$G$-orbits of vertices of $T$   to $\C$: associate to the
orbit of $v\in T$ the equivalence class of
$G_v$. Since any elliptic subgroup is contained in some $G_v$, the
image of $\varphi $ contains $\M$. This shows the first assertion.
Also note that any elliptic subgroup is contained in an elliptic
subgroup whose class is in $\M$. 

The graph
$\Gamma $ has first Betti number   $b_1(\D)$. To prove the
second assertion, it suffices to show that $|\M|$ is an upper bound
for the number of vertices of valence 1 or 2. Suppose $v\in T$
projects onto such a vertex $\ov v$. Since $T$ is minimal and
BF-reduced, the group $G_v$ is a maximal elliptic subgroup, as no
inclusion
$G_e\to G_v$ is onto. Thus $\varphi (G v)\in\M$.  Furthermore,
$\varphi (G v)\neq\varphi (G w)$ if $\ov v$, $\ov w$ are distinct
vertices of valence 1 or 2.  This shows the required bound. 
\cqfd\enddemo

  Assertion (2) of the proposition says that 
{\it accessibility\/}  always holds within a given deformation space
(the only assumption being   $G$   finitely generated). The
accessibility theorem of [\BFa] is equivalent to saying that, given a
finitely presented group $G$, there exists a uniform bound
$|\M|\le C$, valid for all deformation spaces containing a tree with small
edge stabilizers. 

 The trees   such that $s=|\M|$ are precisely the 
  fully reduced trees defined in [\FoGBS]. Forester proved that the
canonical deformation space of a generalized Baumslag-Solitar
group always contains a fully reduced tree [\FoGBS]. He has an example 
showing that this is not true in general. In any case, even when there is
a fully reduced tree, it may not be the most natural element of $\D$.

 \example{Example}
 Let $G=\langle
a_1,b_1,a_2,b_2,a_3,b_3\mid [a_1,b_1]= [a_2,b_2]= [a_3,b_3]\rangle$ 
be the fundamental group of the space obtained by gluing three
once-punctured tori together along their boundary.  Consider the
associated deformation space, where a subgroup is elliptic if and only
if it is contained in a conjugate of $G_i=\langle a_i,b_i\rangle$. The set
$\M$ has 3 elements. There are trees whose quotient graph $\Gamma $ 
has 3 vertices and 2 edges (for instance the tree corresponding to
$G=G_1*_ZG_2*_ZG_3$, with
$Z=\langle[a_i,b_i]\rangle$), but none of them is invariant under the 
automorphism of order 3 mapping $a_i,b_i$ to $a_{i+1},b_{i+1}$ mod 3.
To get an invariant (but not reduced) JSJ   tree, one has to create a
``central'' vertex of
$\Gamma
$ with group
$Z$ (see [\Bo]). 
\endexample

As mentioned above, it is not always true that elements
of
$\M$ are conjugacy classes of maximal elliptic subgroups and
correspond bijectively with the vertices of $\Gamma $ for  
$T$ reduced. Things become nicer   in a \na{} space (see Section
7), or if we work in  a restricted deformation space $\D_\A$  (see
Definition \restri) and consider only  vertex stabilizers not in
$\A$. 

\nom\sommets
\thm{Proposition \sta} For any $T\in \D_\A$, the assignment
$v\mapsto G_v$ induces a bijection from the set of vertices $v$ of   $T$
with
$G_v\notin\A$, to the set of maximal elliptic subgroups not in $\A$. 
\cqfd
\fthm

If we call {\it big\/} a group which is not in $\A$, we get:

\nom\bi
\thm{Corollary \sta} All trees in $\D_\A$ have the same ``big''
vertex stabilizers. \cqfd
\fthm 

\subhead Edge groups\endsubhead

We   associate three families of subgroups of $G$ to a $G$-tree $T$. 
Each family is contained in the next.

$\bullet$ {\it Edge groups\/}: a stabilizer $G_e$, with $e$ an edge.

$\bullet$ {\it Generalized edge groups\/}: a group $H$ such that 
$G_e\inc H\inc G_{e'}$, with $e,e'$ edges of $T$.

$\bullet$ {\it Bi-elliptic groups\/}: a group $H$ contained in some $G_e$
(equivalently, a group fixing two distinct points of $T$).

\example{Example} Consider the trees $T_1,T_2$ in the canonical 
deformation space of
$BS(1,6)$, as in Example \bizz. Fix  $i\in\Z$. The group  $H=\langle
a^i\rangle$   is always bi-elliptic (in both $T_1$ and $T_2$). It is a
generalized edge group if and only if
$i$ divides a power of 6. It is an edge group of $T_1$ if and only if $i$ 
is  a power of 6, an edge group of $T_2$  if and only if $i$ is  twice a
power of 6. This example illustrates the general fact that    bi-elliptic
groups and generalized edge groups are more canonical than edge
groups (see Proposition 4.7).

\endexample

\nom\amin
\thm{Proposition \sta} Let $T$ be reduced. A subgroup $H\inc G$ is
bi-elliptic in $T$ if and only if it is contained in a subgroup  $K$
having one of the
  following forms:
\roster
\item $K=A\cap B$, where $A$ and $B$ are elliptic but $\langle
A,B\rangle$ is not elliptic.
\item $K$ is elliptic and there exists a hyperbolic $g\in G$ such that
$K\inc gKg\mi$. 
\endroster

In particular, all reduced trees in a given deformation space have
the same bi-elliptic subgroups. 
\fthm

\demo{Proof} We first show that any $K$ as indicated  is bi-elliptic.
If
$K$ is as in (1), the groups
$A$ and
$B$ fix distinct points $a,b$, and   $  K$ fixes the segment
$[a,b]$. If
$K$ is as in (2), it fixes a point $x$, and  also  the
point 
  $gx\ne x$. 

Conversely, suppose $H$ fixes an edge $e=ab$. If $G_e$ is properly
contained in both $G_a$ and $G_b$, then (1) holds with 
$K=G_e=G_a\cap G_b$
  (note that $\langle
G_a,G_b\rangle$ is not elliptic). If, say, $G_e=G_a$, there exists a
hyperbolic $g$ such that $b=ga$ (because $T$ is
reduced) and (2) holds with $K=G_a$.
\cqfd\enddemo

\nom\bimorph
\thm{Proposition \sta} Given simplicial (non-metric)   trees $T,T' $ in
the same deformation space, the following conditions are equivalent:
\roster
\item There exist morphisms $f:T\to T'$ and $f':T'\to T$.
\item
$T$ and $T'$ have the same generalized edge groups.
\item
$T$ and $T'$ have the same bi-elliptic groups.
\endroster
These conditions are satisfied if $T$ and $T'$are reduced.
\fthm

See Section 2   for the definition of a morphism.

\demo{Proof} Clearly $(2) \Rightarrow(3)$, as bi-elliptic groups
are just subgroups of generalized edge groups.  It is also clear that $(1)
\Rightarrow(3)$, since a morphism   cannot collapse an edge. Let us
show the  stronger
statement 
$(1)
\Rightarrow(2)$. We have to prove that every edge stabilizer $G_e$ of
$T$ contains an edge stabilizer of
$T'$. Since $f'$ is onto, some edge $e'$ of $T'$ contains a nondegenerate
segment which is mapped injectively into $e$ by $f'$. The  stabilizer of
$e'$ is contained in $G_e$.

We now show $(3)
\Rightarrow(1)$. Given
$T$ and $T'$ in $\D$, there is an equivariant map $f:T\to T'$ which
is linear on edges. It is a morphism if and only if no edge $e$ of $T$ is
mapped to a  point. Suppose $f(e)$ is   a point $v\in T'$.
The group
$G_e$ fixes
$v$, and by (3) it   fixes  a non-degenerate segment
$[v,w]\inc T'$. We redefine $f$ on $e$, by
subdividing $e$ into two halves and mapping each half onto $[v,w]$.
Doing this equivariantly on each orbit of collapsed edges replaces
the original $f$ by a morphism. 

We have seen that (3) holds if $T$ and $T'$ are reduced (Proposition 
\amin). 
\cqfd\enddemo

\definition{Definition \sta{} (bi-morphism class)}
We say that two simplicial trees $T,T'$ are {\it bi-morphically 
equivalent\/}, or belong to the same {\it bi-morphism class\/}, if they
satisfy the equivalent conditions of Proposition \bimorph. Metric trees
are said to be bi-morphically equivalent if the underlying simplicial
trees are. 
\enddefinition

A deformation space $\D$ is thus partitioned
into bi-morphism classes. There is only one class if $\D$ is abelian, as 
any elliptic subgroup is bi-elliptic.

 \nom\aminn
\definition{Definition \sta{} ($\A_{min}$, reduced deformation space 
$\D_r$)} Given $\D$, we denote by $\A_{min}$ the family of subgroups
$H$ described in Proposition \amin{} (the   bi-elliptic groups of any
reduced $T$). They may be characterized as the groups which fix an
edge in every tree belonging to $\D$. The associated restricted
deformation space
$\D_{\A_{min}}\inc\D$, consisting of trees $T\in\D$ with all edge
stabilizers in
$\A_{min}$, is the {\it reduced deformation space\/}, denoted $\D_r$. 
\enddefinition

The space $\D_r$ may be viewed as the bi-morphism class containing 
the reduced trees. Trees in $\D_r$
have the same generalized edge groups, and the same bi-elliptic
groups (those of $\A_{min}$), but not necessarily the same edge
groups. 

Proposition   \findim{} provides a
bound for the complexity of reduced trees. In general, there is no
bound for complexity   in $\D_r$ (it is infinite-dimensional), unless
some descending chain condition holds (for instance, if $\A_{min}$
consists of finite groups). 
  This phenomenon already occurs in the canonical deformation space
of GBS groups, for instance for $BS(2,4)$. There are also  
non-ascending GBS examples (see  the end of
Section 7).

\subhead Peripheral   structure    \endsubhead

Let $\D$ be a deformation space, and $\A_{min}$ as defined above. 
Let $\D_r\inc\D$ be the reduced deformation space,
 so that $T\in\D$ is in 
$\D_r $ if and only if all its edge stabilizers are in  $\A_{min}$.  

We say that a 
vertex stabilizer     is {\it big\/}  if it is not in
$\A_{min}$.  We have
seen (Proposition \sommets{} and Corollary
\bi) that all trees in
$\D_r$ have the same big vertex stabilizers, and that such vertex
stabilizers coincide with maximal elliptic subgroups which are not in
$\A_{min}$.

The goal of this  subsection is to associate a
peripheral structure   $\M_0$ to  a maximal elliptic subgroup $G_0$.  It
depends   only on $\D$ and  
$G_0 $.

First suppose that the maximal elliptic subgroup $G_0$ is big (not
in $\A_{min}$). We consider the set of all
  subgroups $H\inc G_0$ which belong to $ \A_{min}$.
Among such groups, we define   $H\le K$ if $H$ is contained in a
conjugate  
$gKg\mi$ (with $g\in G_0$), and $H\sim K$ if $H\le K$ and $K\le H$.

Let
$\M_0$ be the set of maximal elements of the set of
equivalence classes.  An element
$m\in\M_0$ is a family of subgroups of $G_0$, stable under conjugation
(when no bi-elliptic subgroup of $G_0$ is conjugate to a proper
subgroup of itself, $m$ is just a   conjugacy class of subgroups).

\nom\lesar
\thm{Proposition \sta}
Let $G_0\notin\A_{min}$ be a maximal elliptic subgroup. 
Given $T\in\D_r$, let $v$ be the (unique) vertex with stabilizer $G_0$. 
\roster
\item If an edge  $e$ is adjacent to $v$, there exists $m\in\M_0$
such that $G_e$ is contained in a group belonging to $m$. 
\item Given $m\in\M_0$, there exists an edge $e$ adjacent to $v$
such that $G_e$ belongs to $m$. 
\item The image $\ov v$ of $v$ in the quotient graph $\Gamma =T/G$
has valence at least $|\M_0|$. In particular, $\M_0$ is finite.
\item There exists   $T'\in D_r$ 
such that   equality holds in (3): the valence of $\ov v$ is $|\M_0|$.
\endroster
\fthm

\demo{Proof} 
Let $\C_0$ consist of all  groups $G_f$, for $f$ an
edge of $T$ adjacent to
$v$. It is a finite union of conjugacy classes of subgroups of $G_0$, and
any bi-elliptic subgroup of $G_0$ is contained in an element  of
$\C_0$. Since
$\C_0/{\sim}$ is finite,  
$G_e$ is contained in some
$G_f$ which is maximal in   $\C_0/{\sim}$.
But $G_f$  is also maximal in the set of all bi-elliptic subgroups
of $G_0$, so represents an element   $m\in\M_0$. This proves (1).

To prove (2), let $H\subset G_0$ be a representative of $m\in\M_0$.
Since $H$ is bi-elliptic, $H\subset G_e$ for some edge $e$ adjacent to
$v$. By maximality, $H\sim G_e$.  

Given $m'\neq m$ in
$\M_0$, define $H'$ and $e'$ similarly.  The edges $e$ and $e'$ are in
different  $G_0$-orbits, since otherwise 
$H\sim G_e\sim G_{e'}\sim H'$, contradicting $m\neq m'$. Assertion (3)
follows.

If equality does not hold in (3), there exist oriented edges $e,f$ with
origin $v$, in distinct $G_v$-orbits, with $G_e\inc G_f$. Perform an
elementary expansion, creating an edge $vw$ with stabilizer $G_f$,
with $e$ and $f $ now attached to $w$. The new tree is again in $\D_r$
(no new bi-elliptic group is created), and the valence of $\ov v$ has
decreased by 1. Iterating this operation proves (4). 
\cqfd\enddemo

\definition {Definition \sta{} (peripheral structure $\M_0$)} The set 
$\M_0$ is the {\it peripheral structure\/} of the maximal elliptic
subgroup
$G_0\notin\A_{min}$. 
It describes stabilizers of edges adjacent to
the vertex with stabilizer $G_0$, in any tree belonging to  $\D_r$ (in
particular, in reduced trees). 
  When   no bi-elliptic
subgroup of $G_0$ is conjugate to a proper subgroup of itself, and
equality holds in (3), the set
$\M_0$ is simply the set  of conjugacy classes of   stabilizers of   edges
adjacent to
$v$.
\enddefinition  

 More generally, we now wish to consider any maximal elliptic subgroup
$G_0$, whether in $\A_{min}$ or not.
 The definition of $\M_0$ given above  makes sense if
$G_0\in\A_{min}$, but is not   satisfying since in this case
$\M_0$ consists of a single equivalence class $m=\{G_0\}$. We shall
give a better one.

Let
$G_0$ be
any maximal elliptic subgroup. 
 Given a reduced $T\in\D$, choose $v$
with stabilizer
$G_0$ (it is not necessarily unique). If
$G_e\neq G_v$
 for every adjacent edge $e=vw$,  then $G_0\notin\A_{min}$. Suppose
therefore that there is
 $e=vw$ with $G_e=G_0$. Then $w$ is in the $G$-orbit of $v$ because 
$T$ is reduced,
 and
 $G_w=G_0$ because $G_0$ is maximal elliptic. The projection of   the
fixed subtree $T_0 $  of
$G_0$   into  $\Gamma =T/G$ consists of   loops attached to the 
projection  $\ov v$.   If $N(G_0)$ denotes the normalizer of $G_0$, the
quotient
$N(G_0)/G_0$
  acts on $T_0$ freely, with one orbit of vertices.  It
  is
 free, with   rank $k$ equal to the number of loops.

To define $\M_0$, we consider all
groups $H=G_0\cap G_1$, where $G_1$ is elliptic, not contained in
$G_0$.  Among such
subgroups of
$G_0$, we  define
$H\le K$
 whenever $H$ is contained in   $gKg\mi$ for some element $g $ of the
normalizer
 $N(G_0)$.  The set $\M_0$ is then defined as above, as a set of maximal
elements. This definition agrees with the one given before when
$G_0\notin\A_{min}$.

\thm{Proposition \sta} Let $\D$ be a deformation space, and 
  $G_0 $ be any maximal elliptic subgroup. Let $T\in\D$ be reduced. Let
$v$ be any vertex with stabilizer $G_0$. 
\roster
\item If an edge $e$ is adjacent to $v$, and $G_e\neq G_0$, there exists
$m\in\M_0$ such that $G_e$ is contained in a group belonging to $m$. 
\item Given $m\in\M_0$, there exists an edge $e$ adjacent to $v$
such that $G_e$ belongs to $m$. 
\item The image $\ov v$ of $v$ in the quotient graph $\Gamma =T/G$
has valence at least $|\M_0|+2k$, where $k$ is the rank
 of the free group $N(G_0)/G_0$. 
\endroster
\fthm

We leave the proof to the reader.
 
\head   \sect Topologies on deformation spaces\endhead

In this section, and in the next, $\D$ is a deformation space, or more 
generally a restricted deformation space $\D_\A$.  Because of
Proposition \poin, we assume that
$\D$ is irreducible or genuine abelian.
So far,   $\D$ was mostly viewed as a
set. We shall now view
$\D$ as a topological space and a complex. 

We first consider the set $\T$  of all non-trivial minimal metric 
simplicial
$G$-trees,   up to equivariant isometry.

\subhead{Three       topologies }\endsubhead

The {\it
equivariant  Gromov-Hausdorff topology\/} (or simply Gromov 
topology) on $\T$ is defined as follows. A fundamental system of
neighborhoods for $T\in\T$ is given by   sets
$V_T(X, A,\varepsilon )$, with $X\inc T$  and
$A\inc G$    finite sets, and  $\varepsilon >0$. By definition,  $T'$ is in $
V_T(X, A,\varepsilon )$ if and only if there exists a ``lifting'' map
$x\mapsto \wtilde x $ from
$X$ to $T'$ such that $|d(x,gy)-d(\wtilde x,g\wtilde y)|<\varepsilon $
for every
$x,y\in X$ and $g\in A$.

\nom\passep
\example{ Remark \sta }This topology may fail to be Hausdorff. For 
example, consider
$T$,
$T'$ with the same length function, with $T$   a line  and $T'$   genuine 
abelian. Then every neighborhood of
$T$ contains $T'$. 
\endexample

The {\it axes topology\/} on $\T$ is the coarsest topology making each
translation length function
$T\mapsto \ell_T(g)$ continuous (it is coarser than the
  Gromov  topology, as translation lengths  are continuous
in the Gromov topology).  By [\Pa], the two topologies have the same
restriction to the set of   non-abelian trees (note that dihedral trees
are considered irreducible in [\Pa]). Since non-abelian trees are
determined by their length function, this restriction  is Hausdorff.

We now discuss the weak topology. Beware   that the open or closed
cones to be defined now will be   open or closed in the weak topology,
but not necessarily in the  Gromov or axes topology.
 
The set of metric $G$-trees obtained from a given simplicial  tree $T $
by varying edge-lengths (keeping them positive) will be called an {\it
open cone\/} $\C$. Its dimension $p$ is the number of orbits of edges
(recall that $T$ has no redundant   vertices).  Two metric trees are in
the same open cone if and only if they are equivariantly homeomorphic. 

If we allow
edge-lengths of $T$ to become  $0$ (keeping  at least one of them
positive, so that the tree does not become trivial), we obtain the {\it
closed cone
$\ov\C$\/} spanned   by
$T$.  It consists of trees obtained 
from $T$ by (possibly non-elementary) collapses. It may be 
identified  to an octant
$(\R^+)^p\setminus\{(0,\dots,0)\}$.
  A closed cone  
 is a finite union of open cones.

The {\it weak topology\/} on $\T$   is defined by declaring a set to be
closed when its intersection with every closed cone  is closed. It is 
always finer than the Gromov topology (a weakly converging sequence
converges in the Gromov topology).

If a   closed cone $\ov\C$ meets an abelian deformation space $\D$,
then it is contained in $\D$, as any collapse or expansion 
 of an abelian tree is  
elementary (in the sense of Section 3).  
If $\D$ is irreducible, then $\D\cap\ov\C$ is
not necessarily closed in $\ov\C$ (this will be discussed
later), but it is convex. To see this, we use Condition (7) of Theorem
\lessix :  if $\ell$ and $\ell'$ are bi-Lipschitz equivalent, then any
convex combination is also equivalent to $\ell$ and $\ell'$.

\subhead Projectivized space  \endsubhead

Rescaling the metric defines an action of the group of positive reals 
on $\T$, and the real object of interest is the quotient space $P\T$. It is
a complex, an open   (resp\.   closed) simplex  being the projection of
an open (resp\. closed) cone of $\T$. The three
topologies defined on
$\T$ induce quotient topologies on
$P\T$. 

If we choose  a finite generating system $S$ for $G$, 
the projection $\T\to P\T$
has a natural  section, defined as the set of trees such that
$\ds\sum_{s\in S}\ell(s)+\sum_{s,t\in S}\ell(st)=1$ (recall that a tree
such that all elements 
$s$ and $st$ are elliptic is trivial [\Ser]).   This section is continuous in
all three topologies (because the maps
$T\mapsto\ell_T(g)$ are), so
$\T$ is homeomorphic to
$P\T\times\R$. In practice, we will work in either 
$\T$ or $P\T$, and the results will automatically apply to the other 
space. 

\subhead{Topologies on $\D$}\endsubhead

We now   restrict  these   topologies to a given deformation space
$\D$ (irreducible or genuine abelian). 
First consider the Gromov topology and the axes topology.

If $\D $ is
irreducible, the restrictions are the same. The induced topology  will 
be called the
  Gromov topology on $\D$.  It is Hausdorff.

If $\D$ is genuine abelian, all trees in $\D$ have the same length 
function (up to scaling) and the axes topology is trivial on $P\D$. We
therefore only consider the   Gromov  topology. We shall prove that it
is Hausdorff (it is not obvious in this case), see Proposition 5.7. 

We now consider the weak topology. There is a subtlety here in the  
irreducible case, due to the fact that $\D$, though a union of open cones
of
$\T$, is not necessarily a union of closed cones ($\D$ is not closed in
$\T$ for the weak topology).

Indeed, the closed cone
$\ov\C$ of
$\T$ spanned   by a given  
$T\in\D$ is not always entirely contained in $\D$, as collapsing certain 
orbits of edges may yield a tree not in $\D$.
 The intersection of $\ov\C$ with
$\D$ is the set of trees that may be obtained  from trees in $\C$ by 
elementary collapses (as defined in Section 3). It will be called a {\it
closed cone of\/}  
$\D$ (it may have faces ``at infinity''). Closed simplices of
$P\D$ are defined similarly. 
A closed cone (resp\. closed simplex) is a finite union
of open cones (resp\. simplices). It always contains reduced trees.

 \example{Remark}   Although simplices of $P\D$ are not actual
simplices, there is a standard way to obtain   a genuine simplicial
complex (see [\CVa,
\McM]). Start from the barycentric subdivision of $P\T$, and look at
the union of closed simplices which are contained in $P\D$. One easily
checks that one obtains a simplicial complex contained in $P\D$ on
which
$P\D$ retracts by deformation.
\endexample

The weak topology on $\D$ is defined as the induced topology. 
A subset of $\D$   is closed in $\D$  if and only if its intersection with 
any closed cone of $\D$ is closed (in the cone).  Similarly, a subset of
$P\D$   is closed  in the weak topology if and only if its intersection
with any closed simplex  of $P\D$  is closed.

As pointed out in [\McM], the weak topology does not always coincide
with the Gromov topology or the axes topology. Indeed, the weak
topology on a locally infinite complex is never metrizable (there is no
countable basis of neighborhoods).

\nom\memetop
\thm{Proposition \sta} Let $\D$ be a deformation space. The Gromov
topology and the weak topology induce the same topology on any
finite union of cones of $\D$ (resp\. of simplices of $P\D$).
\fthm

The proof   will show that the analogous result for cones in  $\T$ is
true if none of the cones is contained in an abelian deformation space.
But it may be false on the
union of two cones contained in different abelian deformation spaces
(see Remark
\passep). The proof of the proposition for $\D$ abelian requires the
fact that $\D$ is Hausdorff in the Gromov topology. This will be
   proved at the end of this section.  

\demo{Proof} It suffices to show the result in $P\D$.  Let $\Sigma
$  be any finite union of simplices of  $P\T$ (not necessarily contained 
in $P\D$). Let
$\ov\Sigma
$  be  its weak closure in
$P\T$. Trees in
$\ov\Sigma $ are obtained by performing (possibly
non-elementary)   collapses on trees in
$\Sigma $. Being a finite  union of closed simplices of
$P\T$,  the space
$\ov\Sigma
$ is weakly compact. Since the weak topology is finer than the Gromov
topology, it is enough to show that $\ov\Sigma $ is Hausdorff in the 
Gromov topology. 
 
If $\Sigma$ contains only non-abelian trees, so does   $\ov\Sigma $,
  and we know that the set of non-abelian trees is
Hausdorff in the Gromov topology. If $\Sigma $ is contained in a genuine abelian
deformation space $\D$,  
we have seen   that  $\ov\Sigma =\Sigma $. We
then use the fact that $P\D$ is Hausdorff in the
Gromov topology 
(see Proposition 5.7).
\cqfd\enddemo

\nom\bonnetop
\example{Remark \sta}  In particular, the topologies agree on  any
open or closed cone of $\D$.   Such a cone being convex (as a subset of
$ \R ^p$), it is contractible.  The same holds for simplices of $P\D$.
\endexample

\subhead{Locally finite   trees}\endsubhead  

In this subsection we assume that all trees in $\D$ are locally finite
(see  Section 4 for a discussion,    and examples of local finiteness). 
In this case, it is easy to show that the complex $P\D$
is locally finite (see [\LeGBS]). 
Moreover,   all vertex and edge stabilizers
are   commensurable, so they are all finitely generated, or all  
  infinitely generated. We shall
  show:

\nom\TopoLocFin
\thm{Proposition \sta{}} If $\D$ consists of locally finite trees with 
finitely generated  vertex stabilizers, then the
Gromov topology and the weak topology coincide on $\D$ and $P\D$. 
\fthm

Let $T$ be a   tree.
Denote by $E(T)$ the set of  oriented open   edges of
$T$, and by
$e\mapsto
\ov e$ the involution mapping an edge to the opposite  edge. Say
that three oriented edges   $e_1,e_2,e_3$ are aligned in this order
 if there is a geodesic going through $e_1,e_2$ and $e_3$
successively, with the correct  orientation.

 \nom\collap
\thm{Lemma \sta{}} Let $T,T'$ be two simplicial $G$-trees. Then $T$ is a 
(maybe non-elementary) collapse of $T'$ if and only if there exists a 
$G$-equivariant 
    injection $\phi:E(T)\hookrightarrow E(T')$, with 
$\phi(\bar e)=\bar{\phi(e)}$,  which preserves alignment in the
following sense: if $e_1,e_2,e_3$ are aligned in this order, then so are
$\phi(e_1),\phi(e_2),\phi(e_3)$.
\fthm

\demo{Proof}  If $f:T'\ra T$ is a collapse map, then the preimage of each
open   edge of $T$ is an open edge of $T'$. This defines a map
$\phi:E(T)\ra E(T')$ with the desired properties.

Conversely, we are given $\varphi $   and we want to define a
collapse map $f:T'\ra T$. Let
$I\inc T'$ be the union of all open edges in the image of $\varphi $. We
first define
$f$ on
$I$,  using
$\phi\m$. Next,   consider $x\in T'\setminus I$. Let
$\alpha,\beta$ be two edges in $I$ such that
$x$ is in the convex hull of $\alpha\cup \beta$ (they   exist because
$T'$ is minimal). Choose
$\alpha,\beta$ so that this hull is minimal. Orient $\alpha$ and $\beta$
so that they   point towards  $x$. Denote by $y_\alpha$ and
$y_\beta$ the terminal points  of $\phi\m(\alpha)$ and $\phi\m(\beta)$
respectively. If
$y_\alpha\neq y_\beta$,   there is an edge $e $
between them, and  $e\neq \alpha ,\beta$  because $\varphi  $
preserves alignment.
 Then
$\phi(e)$ is between $\alpha$ and $\beta$,   contradicting minimality
  in the choice of $\alpha,\beta$. Thus
$y_\alpha= y_\beta$, and we may  define 
$f(x)=y_\alpha$ consistently. The map $f$ is clearly  a collapse map.
\cqfd\enddemo

\demo{Proof of   Proposition \TopoLocFin}   
If $T, T'$ are metric trees, we say that $T'$ collapses to $T$ if this is
true for the underlying simplicial trees.  Local finiteness implies
that, given
$T\in\D$, the set $S(T)$ consisting of  trees of $\D$ which collapse to
$T$ is   a finite union of cones.
  Since the weak topology is finer than the Gromov topology, and
both topologies coincide on $S(T)$ by Proposition \memetop, it
 suffices to prove that $S(T)$ contains the intersection of $\D$ with
some  Gromov-neighborhood
$V_T(X,A,\eps)$ of $T$.

By
local finiteness of $T$, one can find  a finite subtree $K$ of $T$ such
that,  for each pair $\alpha,\beta$ of adjacent edges, there exists $g\in
G$ mapping both $\alpha$ and $\beta$ into $K$. Let $X$ be the set of
vertices of $K$. Let $A\inc G$  be a finite set containing a set of
generators of each $G_v$ for $v\in X$ (recall that each $G_v$ is
assumed to be finitely generated), and also an element sending
$v_1$ to
$v_2$ for each pair of points $v_1,v_2\in X$ lying in the same
$G$-orbit. Take
$\eps$ very small compared to the lengths of edges of $T$. We show
that any $T'\in \D\cap V_T(X,A,\eps)$ collapses to $T$.

 We first construct an equivariant map
$v\mapsto  v'$ from $V(T)$ (the set of vertices of $T$) to $T'$. If  $v\in
X$, we denote by $\wtilde v$ a ``lift'' of $v$ to
$T'$, as in the definition of  $ V_T(X,A,\eps)$. For
each $G$-orbit of vertices, choose a representative $v_i\in X$. 
Since $G_{v_i}$
is elliptic in
$T'$, and   the generators of $G_{v_i}$ move $\wtilde v_i$ by less
than
$\eps$, the group 
$G_{v_i}$ has a fixed point in $T'$ at distance at most $\eps/2$ from
$\wtilde v_i$. We define $v'_i$ to be such a point, and we extend
equivariantly.

If $v\in X$ is not a $v_i$, the point $v'$  is at distance at most
$3\eps/2$ of $\Tilde v$. It follows that, if $uv$ and $uw$ are two adjacent
edges of
$T$, then  the distances between $  u',  v', w'$ are
$4\eps$-close to the distances between $u,v,w$ (this is true if the
points are vertices  of  $K$, hence in general by  our choice of $K$).
In
particular, the overlap between $  u'  v'$ and $u'
w'$ has length at most $6\eps$.

We shall now define a map $\phi:E(T)\ra E(T')$ as in Lemma \collap. For
each edge
$e=uv$ of
$T$, let
$I_e$ be the segment of $T'$ obtained by removing from $
u' v'$ the $ 6\eps$-balls   centered at $u'$ and $v'$. Our choice of
$\eps$ guarantees that  $I_e$ has length much bigger than
$\eps$. If $e_1,e_2$ are adjacent, then $I_{e_1}\cap I_{e_2}$ is empty
and 
$d(I_{e_1},I_{e_2})\le12\eps$.  
It follows that, given a geodesic edge path
$e_1,\dots,e_n$ in $T$, the segments
$I_{e_1},\dots,I_{e_n}$ are aligned in this order in $T'$. In particular,
$I_{e_1}\cap I_{e_2}=\es$ for any pair of edges $e_1\neq e_2$,.

Moreover,   no $I_{e}$ contains a vertex in its interior. Otherwise,
the subtree  branching from $I_e$ would contain  a
segment $I_{e'}$ (by minimality). We would get a contradiction by
considering the edge path from $e$ to $e'$ in $T$, and the
corresponding aligned segments in $T'$. 

Thus, one can define
$\phi:E(T)\hookrightarrow E(T')$ by mapping an edge $e$ to the edge of
$T'$ containing $I_e$. By the lemma, $T$ is a collapse of $T'$.
\cqfd\enddemo 

\subhead{The closure  of an irreducible space}\endsubhead

Let $\D$ be an irreducible deformation space. 
We   consider its closure $\ov
\D$   in the space $\T_\R$ of non-trivial minimal \Rt s, equipped with
the Gromov topology. It usually contains non-simplicial \Rt s  (see
[\CVb] for a concrete example). 

If $T\in\ov\D$, we denote by $\ell$ its length
function. 
  Given $g\in G$, we
denote by $A_g$ the   axis of $g$ if $g$ is hyperbolic in $T$, its fixed
point set if it is elliptic (an element which is hyperbolic in $\D$ may
become elliptic in $T\in\ov\D$).  

\nom\obord
\thm{Proposition \sta} Let $\D$ be an irreducible deformation space.
Let $\ov\D$ be its closure in the Gromov topology. Let $T$ be an
\Rt{} in
$\ov\D$. 
\roster
\item If $g\in G$ is hyperbolic in $\D$ but elliptic in $T$, its
fixed point set contains no tripod.
\item If $g,h,[g,h]$ are hyperbolic in $\D$, then   $A_g\cap
A_h$ is a segment (possibly empty or degenerate) whose length is at
most
$\ell(g)+\ell(h)$. 
\item   $T $ is irreducible. 
\item $\ov\D$ is also the closure of $\D$ in the axes topology.
\endroster
\fthm

\demo{Proof}
Assertion (1) holds because  fixing a tripod is an open condition 
(see   [\PaENS, p\. 153]). To prove (2), note that $A_g$ and $A_h$
contain no tripods.   If
$|A_g\cap A_h|>\ell(g)+\ell(h)$, this inequality holds in trees $T'\in\D$
close enough to
$T$ (if $x\in A_g$, any lift $\wtilde x$ of $x$ to $T'$ is  close to the
axis of
$g$ in
$T'$). 
A standard argument shows that $[g,h]$ is elliptic in $T'$, a
contradiction. 

Since $\D$ is irreducible, we can find $g,h$ as in (2). 
Compactness of $A_g\cap A_h$  implies that $T$ is not abelian.
Applying the argument to $g^2$ and $h^2$ shows that $T$ is not
dihedral, so it  is irreducible.

   Recall that the Gromov topology and the axes  topology agree
on the set of non-abelian trees [\Pa]. To prove  (4),   it suffices to
see that length functions of trees in
$\D$ cannot accumulate onto an abelian length function $|\varphi |$, with
$\varphi :G\to\R$ a nontrivial homomorphism.  

Supposing they do,   fix
$g\in G$ with
$\varphi (g)>0$. It is hyperbolic in $\D$. Since
$\D$ is irreducible,  we can find a conjugate $h$ of $g$ such that $[g,h]$
is hyperbolic in $\D$ (choose $h$ such that $g$ and $h $ have disjoint
axes in some tree of $\D$). Now consider $T\in\D$ whose length
function $\ell$ is close to $|\varphi |$ when applied to $g$, $h$,
 and  
$g^3h^{-3}$ (in particular,
$\ell (g^3h^{-3}) $ is small). Using the formulas in [\Pa, Proposition
1.6], one sees that the intersection of the axes of
$g$ and
$h$ in $T$  must contain a segment of length almost $ 3\varphi (g)$.
Since
$3\varphi (g)>2\varphi (g)\sim\ell(g)+\ell(h)$, the commutator 
$[g,h]$ is elliptic in
$T$, a contradiction. 
 \cqfd\enddemo

\subhead{Abelian deformation spaces   are Hausdorff}\endsubhead

\nom\AbelianHausdorff
\thm{Proposition \sta}
Any genuine abelian deformation space  $\D$ is Hausdorff in the Gromov
topology.
\fthm

\demo{Proof} 
If $T\in\D$, all translation axes $A_g$ in $T$ have a
common end, which is fixed by $G$.  We say that  a hyperbolic $g$ is
positive (in $T$) if it moves points away from the fixed end, negative if
it moves towards the end (in other words, we think of the fixed end as
being at $-\infty$). 

It is easy to check that a  hyperbolic element $g$ is positive if and only
if, given any
$h\in G$, the group generated by the commutators   
$[g^{n},h]$,  
   with $n>0$, is elliptic. In particular, positivity
does not depend on the choice of $T$ in $\D$. 

Similarly,   the fact that two  hyperbolic elements $g,h$ have the same
axis depends only on $\D$   (it is characterized by ellipticity of the
group generated by the elements   $[g^{ n},h]$, $n\in\Z$).

Consider positive hyperbolic elements $g,h$ with different axes. Let
$P_{gh}$ be the endpoint of the ray $A_g\cap A_h$   (see  Figure 1). It
is characterized by the equations 
$$ \left\{\aligned d(P,gP)&=\ell(g)  \\ d(P,hP)&= \ell(h)\\
d(gP,hP)&=\ell(g)+\ell(h).
\endaligned  \right.\tag{$1_{gh}$ }$$ 

\midinsert
\centerline{\includegraphics[scale=1]{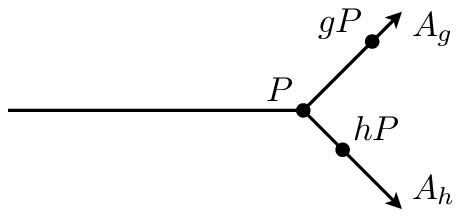}} 
\botcaption{Figure 1}{Characterizing    the endpoint of $A_g\cap A_h$ }
\endcaption  
\endinsert

Furthermore, this system is stable in the following sense: there exists
a universal constant
$C$ such that, if
$P$ satisfies each  equation  up to
$\varepsilon $, then it is $C\varepsilon $-close to $P_{gh}$.

Now consider   three   positive hyperbolic elements  $g,h,i$ with
distinct axes. We claim that the number $\delta _T =d(P_{gh},P_{gi})$
depends continuously on
$T$ (in the Gromov topology).  

 If $T'$ is another tree in $\D$, with
length function $\ell'$, we denote by
$(1'_{gh})$ the system $(1_{gh})$ with $\ell $ replaced by $\ell'$. Note
that the right-hand sides   depend continuously on $T$. 
Fix
$\varepsilon >0$. If
$T'$ is close enough to
$T$, we can find    lifts  $\wtilde P_{gh}\in T'$ satisfying $(1'_{gh})$
up to
$\varepsilon $, and 
$\wtilde P_{gi} $ satisfying $(1'_{gi})$ up to $\varepsilon $, with
$d(\wtilde P_{gh}, \wtilde P_{gi})$ $\varepsilon $-close to $\delta _T$.
We get
$|\delta _{T'}-\delta _T|\le(2C+1)\varepsilon $, showing continuity.

To prove that $\D$ is Hausdorff, suppose that $T ,T'$ do not
have disjoint neighborhoods. Then $ \delta _{T'}=\delta _T$ for all
ordered triples $(g,h,i)$ as above. We show that this implies $T =T'$. 

Denote by $g_1,g_2,\dots$ the positive hyperbolic elements of $\D$.
Write
$A_i$ (resp\. $A'_i$) for the axis of $g_i$ in $T$ (resp\. $T'$).  We may
assume
$A_1\neq A_2$.   Since
$A_i=A_j$ is equivalent to $A'_i=A'_j$, we get $A'_1\neq A'_2$. 

There is a unique isometry $f_2:A_1\cup A_2\to A'_1\cup A'_2$
sending $A_1$ to $A'_1$ and $A_2$ to $A'_2$. If $A_3$ is distinct from
$A_1$ and $A_2$, its position   with respect to $A_1\cup A_2$ is
completely determined by the three numbers $\delta  _T$ associated
to the triples
$(g_1,g_2,g_3)$, $(g_2,g_3,g_1)$, $(g_3,g_1,g_2)$ (one of these
numbers is $0$, the other two are equal). Thus $f_2$ extends 
uniquely to an isometry $f_3:A_1\cup A_2\cup A_3\to A'_1\cup
A'_2\cup A'_3$ sending $A_3$ to $A'_3$. 

By a similar argument,  $f_3$ has a sequence of successive unique
extensions $f_n: A_1\cup\dots\cup A_n\to A'_1\cup\dots\cup A'_n$
sending $A_i$ to $A'_i$ for $i\le n$, and  we finally    get a global
isometry
$f:T\to T'$ sending every $A_g$ to $A'_g$.   By uniqueness,
$f$ is $G$-equivariant.  

\cqfd\enddemo

\head \sect Contractibility \endhead

Let $\D$ be a deformation space, or more generally a restricted 
deformation space
$\D_\A$. We shall prove:

\nom\contract
\thm{Theorem \sta} Let $\D $ be a  deformation space.
\roster
\item $\D$ is contractible in the weak topology.
\item If $\D$ contains  a tree
$T_0$ with finitely generated vertex stabilizers, then $\D $ 
is contractible in the  Gromov topology.
\item If $\D$ is irreducible and contains  a tree
$T_0$ with finitely generated vertex stabilizers, then  its
closure $\ov\D$  in the  
Gromov  topology is  contractible.
\endroster The same results hold for $P\D$. 
\fthm

Recall that $\ov\D$ is the closure of $\D$ in $\T_\R$ (the space of
non-trivial minimal \Rt s, equipped with the Gromov topology or the
 axes topology, see Proposition \obord). Its projectivization is the
closure of
$P\D$ in
$P\T_\R$.  Note that the closure in $\T$ (equipped  with the Gromov
topology) of a genuine abelian space   is never Hausdorff, as adding  
the linear abelian tree having the same length function makes
$\D$ non-Hausdorff (see Remark \passep).

The key technique for proving contractibility is Skora's idea of
deforming morphisms between metric trees: given any morphism
$f:S_0\to S
$, there is a canonical way of constructing {\it intermediate trees\/}
$S_t$, for $0\le t\le \infty$, with
$S_\infty=S$    (see [\Cl, \GL, \Sk]).    The tree
$S_t$ depends continuously on
$f$ and $t$ in the Gromov topology. It belongs to $\D$ (resp\. $\D_\A$)
if both $S_0$ and
$S$ do, because there are morphisms $\varphi _t:S_0\to S_t$ and $\psi
_t:S_t\to S$ with $\psi
_t \circ\varphi _t=f$. 

\subhead{Defining a contraction}\endsubhead

In this subsection, we explain how to define a map $\rho
:\D\times[0,+\infty]\to\D$ with   $\rho
(T,\infty)= T$  and $\rho (\D\times\{0\})$ contained in a contractible
subset. We will then show Assertions (1) and (2) by establishing
continuity of $\rho $.  The proof of (3) will be
given at the end of the section.

Fix $T_0\in\D$, with finitely generated vertex stabilizers if proving
(2). We may assume that it
is reduced (elementary collapses do not create  new vertex
stabilizers). We view $T_0$ as a simplicial tree, and we let $\C_0\inc\D$
be the open cone containing it (it is the set of all metric trees with
$T_0$ as underlying simplicial tree). Recall (Remark \bonnetop) that
$\C_0$ is contractible, both in the Gromov topology  and in the weak
topology.

We shall associate to $T\in \D$  a
morphism
$f_T:T_0(T)\to T$, where $T_0(T) $ is a metric
tree belonging to $\C_0$. 
Skora's deformation   (see  above)  
provides intermediate trees $T_t(T)$. Setting $\rho (T,t)=T_t(T)$ 
yields a   map
$\rho :\D\times[0,+\infty]\to\D$, with $\rho (T,0)=T_0(T)$ and $\rho
(T,\infty)= T$   (if some
intermediate tree fails to be minimal, we replace it by its minimal
subtree). 
Since $T_0(T)$ belongs to the contractible set $\C_0$, we deduce that
$\D$ is contractible  (assuming continuity of $\rho $).

There remains   to construct $f_T$.
It  is   easy in the following special case (see [\GL]). 
Suppose that  $\D$ is a restricted deformation space
$\D_\A$, and   no vertex stabilizer
$G_v$  of $T_0$ belongs to $\A$. Then $G_v$ fixes a {\it unique\/} point
$v'$ in $T$, and
$f_T$ is defined on the vertex set of $T_0$ by sending $v$ to $v'$.  
Since
$T_0$ is reduced, adjacent vertices are mapped to distinct points by
Remark \pasecr. There is a unique edge-length assignment on $T_0$
such that the natural extension of $f_T$ is an isometry when restricted
to edges. 
This assignment defines the  metric tree
$T_0(T)\in\C_0$, and $f_T:T_0(T)\to T$ is a morphism. 

In general (as in Skora's proof that Culler-Vogtmann space is
contractible), constructing $f_T$ requires     choosing  a
basepoint
$P\in T$ (continuously). Once 
$P$ has been chosen, one obtains $f_T$ as follows. Fix a
representative 
$v_j\in T_0$ in each orbit of vertices and define $f_T(v_j)\in T$ as the 
projection
$P_j$ of the basepoint $P$ onto the fixed point set of
$G_{v_j}$. Then extend $f_T$ equivariantly to the vertex set of $T_0$
 (this is possible because $G_{v_j}$ fixes
$P_j$), and to the edges as before.

Skora's
construction of a basepoint was based on minimization. We give a
different one.

 First suppose that   $\D$ is
irreducible. Fix elements $g,h\in G$ such that
$g$, $h$, and $[g,h]$ are hyperbolic in $\D$. If $T\in\D$, the axes $A_g$
and $A_h$ of $g$ and $h$ in $T$ have   compact intersection (see
Proposition \obord), and we use this fact to define $P$. 

The most symmetric way would be to  define
$P$ as the midpoint of the segment $A_g\cap A_h$ if $A_g$ and $A_h$
meet, as the midpoint of the ``bridge'' between them if they  are
disjoint. This ``symmetric basepoint'' will be used at the end of this
section (see Figure 4), but  for technical reasons (in the proof of
Corollary 6.3) we now define
$P$ in  a slightly different way (see Figure 2 later in this section). If
$A_g$ and
$A_h$ meet, we order Ag so that the action of g is by a positive
translation and we let P be the largest element of the segment
$A_g\cap A_h$. If
$A_g$ and $A_h$ are disjoint, we define $P$ as the point of $A_g$
closest to $A_h$. 

If $\D$ is abelian, 
we fix two positive hyperbolic elements $g,h\in G$ with
distinct axes in some and therefore in any $T\in \D$ (see the proof of
Proposition 
\AbelianHausdorff). We define
$P$ as the endpoint of the ray $A_g\cup A_h$.

This completes the definition of the map $\rho
:\D\times[0,+\infty]\to\D$. 
We now study its topological properties. 

\nom\critere
\thm{Lemma \sta}  Suppose that, given any  two vertices $v,w$ of
$T_0$, the map
$T\mapsto d(f_T(v),f_T(w))$, from $\D$ equipped with the Gromov
topology to
$\R$, is continuous. Then the  map
$\rho $ is continuous in the Gromov topology.
\fthm

Recall that all trees $T_0(T)$ have $T_0$ as their underlying simplicial
tree, so their vertex sets are canonically identified. 

\demo{Proof} If the condition of the lemma is satisfied, it is easy to
see that the morphism $f_T$ depends continuously on $T$ (there is a
natural Gromov topology on the space of morphisms, see e.g\. [\GL]).
Continuity of
$\rho $ then follows from continuity of  Skora's deformation (see
[\GL, Section 3]). We may have to replace a tree by its minimal
subtree, but this  is a continuous operation in the Gromov topology. 
\cqfd\enddemo

\nom\cg
\thm{Corollary \sta} If $\S$ is a closed cone of $\D$, the restriction
$\rho _\S:\S\times [0,\infty]\to\D$ is continuous in the Gromov
topology.
\fthm 

\demo{Proof}  By Proposition \memetop, the Gromov topology
and  the weak topology   agree on
$\S$. It is homeomorphic to a Euclidean cone, parametrized by
edge-lengths $\ell_1,\dots,\ell_p$ ($p$ is the number of orbits of
edges of trees in the interior of $\S$; as discussed in Section 5, the
cone may fail to be  closed as a Euclidean cone). The basepoint
$P$, and its projections $P_j$, were defined using only the simplicial
structure of
$T$, with no reference to edge-lengths. In restriction to $\S$, any
function
$d(f_T(v),f_T(w))$ as in Lemma \critere{} is therefore a fixed linear
combination of $\ell_1,\dots,\ell_p$, hence is continuous.
\cqfd\enddemo

\subhead Contractibility in the weak topology\endsubhead

We now show that $\rho $ is continuous in the weak topology, thus  
proving contractibility of $\D$ in that topology. By definition of the
weak topology, it suffices to  prove continuity of restrictions 
$\rho _\S:\S\times [0,\infty]\to\D$, with the target space equipped
with the weak topology. 

  As in [\GL], we show
this by proving that, when
$T$
 varies inside $\S$, the set of intermediate trees $T_t(T)=\rho (T,t)$
only meets finitely many cones. 

\nom\finit
\thm{Lemma \sta{} (finiteness lemma)}
Let $\C$ be an open cone in $\D$. The set of intermediate trees
$T_t(T)=\rho (T,t)$, for $t\ge0$ and $T\in\C$, is contained in a finite
union of cones. 
\fthm

The meaning of this lemma is that  there are only finitely many
possibilities for
$T_t(T)$ as a non-metric tree, as $t$ and $T$ vary. The lemma has
weaker hypotheses than Lemma 4.3 of [\GL], and the conclusion is
weaker: we do not control all trees $T'$ such that $f_T$ factors through
$T'$, only those which occur as intermediate trees in the deformation. 

Since the closed cone $\S$ is a finite union of open cones,
contractibility of $\D$  in the weak topology follows from the finiteness
lemma, together with Corollary \cg{} and Proposition  \memetop. 

The proof of Lemma  \finit{} is rather complicated, so  we first  
prove the following simpler result:

\nom\finitd
\thm{Lemma \sta}
Let $f:T_0\to T$ be a morphism between two metric trees belonging
to the same   deformation space. The set of  intermediate trees 
$T_t$ is contained in a finite union of cones.  
\fthm

\demo{Proof} 
The proof we   give is not the shortest possible, but the method  will be used
to prove Lemma \finit. We first recall how  the intermediate
trees $T_t $, and the factorization of  $f $   through
morphisms    $\varphi _t:T_0 \to T_t $ and $\psi _t:T_t \to T $, are
constructed
  (see  [\GL, \Sk]). 

Two points $x,y\in T_0 $ with
$f (x)=f (y)$  are identified in
$T_t $ if and only if the image of the segment $[x,y]$ by $f $ is contained in the
$t$-ball centered at $f (x)$. This defines $T_t$ as a set. The distance 
function on
$T_t$ is the maximal one making the quotient map $\varphi _t$ 
$1$-Lipschitz, and
$\psi _t$ is defined in the obvious way. See [\GL] for details. 

A key feature of this
construction is that the geometry of the image of a segment $[x,y]$ in
$T_t$ depends only on $t$ and the restriction of $f$ to $[x,y]$. 

By subdividing $T_0$ and $T$, we may assume that $f$ is a simplicial 
map (an edge is mapped to an edge). 
Choose a set of representatives $\varepsilon _i$ for orbits of
 non-oriented edges of $T$. Label edges of $T$ by the
corresponding $i$. Also orient edges of $T$ in an equivariant way. Use
$f$ to lift labeling  and orientation to edges of $T_0$.

Consider a couple   $(\alpha,\beta)$ consisting of distinct
non-oriented   edges of
$T_0$ with the same image in $T$. Let   $E$ be the convex hull of
$\alpha\cup \beta$ (a segment). Since $T_0$ and $T$ have the same
elliptic subgroups,   Remark
\bdd{}   implies that   there is a   bound for the   length of $E$ (as
$(\alpha ,\beta )$  varies).

The morphism $\varphi _t:T_0\to T_t$ folds $E$ into  a finite subtree 
$E(t)$ which interpolates between $E$ and its image
$F=f(E)$  in $ T$.   As already pointed out, this folding process
depends only on the restriction of
$f $ to $E$.  
We formalize this observation as follows.  
  
 Let $(\alpha,\beta)$ and $(\alpha',\beta')$ be two couples as above.  
We say that they are {\it equivalent\/} if there exists a commutative
diagram 
$$\xymatrix{E\ar[r]\ar[d]_f&E'\ar[d]^f\\ F\ar[r]&F' }$$ where the
horizontal arrows are label-preserving and orientation-preserving  
simplicial isomorphisms.  In particular,  $(g.\alpha,g.\beta)$ is
equivalent to $(\alpha,\beta)$ for any $g\in G$. Also note that  the
common label of $\alpha$ and $\beta$ gives a well-defined label $i$
   for the corresponding equivalence class.

Equivalence implies that  the same folding takes place in $E$ as in $E'$ 
as $t$ increases; in particular, $E(t)$ and $E'(t)$ are isometric. 
Furthermore,   there are only finitely many equivalence classes since
$E$ has bounded length and there are finitely many labels.

  Given $t>0$ and an equivalence class $c$, we
shall define a 
(possibly empty or degenerate) subsegment  $I_c(t)$ of the edge 
$\varepsilon _i$ (where
$i$ is the label of $c$). It will encode the way in
which  folding occurs in 
$T_t$. 

Represent $c$ by $(\alpha,\beta)$, with   $\alpha,\beta$  mapping onto
$\varepsilon _i$. The edges $\alpha,\beta$  
are mapped isometrically into
$T_t $, but their images may be glued along a closed segment 
(possibly a point). We define $I_c( t)$ as  the image of this segment in
$\varepsilon _i$ (it is empty if no gluing occurs). It depends only on
$c$, not on the choice of
$(\alpha,\beta)$.

We have associated to $t$ a finite family of subsegments $I_c( t)$,  one
for each equivalence class $c$.   We  focus on the combinatorial
structure of this family: we say that $t$ and $t'$ are similar if, for every
$i$, there exists an orientation-preserving homeomorphism of
$\varepsilon _i$ mapping each 
$I_c( t)$ to $I_c( t')$. There are only finitely many similarity classes. We
complete the proof by fixing a   class and showing that, as a non-metric
tree,
$T_t $  does not depend on   the choice of $ t $ in that class.

  Subdivide each $\varepsilon _i$ by adding all endpoints of
segments $I_c$. Extend the subdivision equivariantly, and lift it  to 
$T_0$ and $T_t$ using $f$ and $\psi _t$. We use the symbol $\
\widetilde{}\ $ to indicate the subdivided trees. The simplicial
structure on   $\widetilde T_0$ is fixed, and the map
$\varphi _t$ is now simplicial (it sends edges to edges). The simplicial structure on
$\widetilde T_t $ is completely determined by knowing which pairs of edges of
$\widetilde T_0$ get identified in
$\widetilde T_t $ (see [\GL]). 

So consider edges $\widetilde \alpha, \widetilde \beta$
with the same image  $\widetilde\varepsilon $ in $\widetilde T$. They 
are contained in edges
$\alpha,\beta$ of
$T_0$ (the unsubdivided tree). Using $G$-equivariance, we may
assume that
$\alpha$ and
$\beta$ map onto some $\varepsilon _i$. Now observe that  $\widetilde
\alpha$ and $
\widetilde \beta$ are identified in $\widetilde T_t $ if and only if
$\widetilde\varepsilon $ is contained in $I_c( t)$, where $c$ is the 
equivalence class of $(\alpha,\beta)$. This description  does not
involve $t$ (in the given similarity class), so $\widetilde T_t $ depends
only on the similarity class (as a   simplicial tree).
\cqfd\enddemo

\demo{Proof of Lemma \finit}  Choosing a point in  $\C$ amounts to
choosing a metric
$m$ (which assigns lengths to edges) on a fixed non-metric tree. 
Instead of controlling the image of a single curve $t\mapsto T_t$ as in
the previous lemma, we now have to control a family of curves
$t\mapsto T_t(m)$, indexed by $m$. We review the arguments in the
proof of Lemma \finitd, checking that the   parameter $m$ introduces
only finitely many new possibilities. 

In order to make $f_T$ simplicial, we only
  have to subdivide   $T_0(T)$.  Because of the way $f_T$ was
constructed (using the basepoint $P$ and its projections   $P_j$),
the
subdivision of $T_0(T)$, and $f_T$ as a simplicial map,  are  
independent of $m$. In particular,  equivalence classes of couples
$(\alpha,\beta)$ are defined independently of $m$. The subsegments
$I_c(t)$ now depend on $m$, so we denote them by 
$I_c(t,m)$. Similarity is   defined as in the proof of Lemma \finitd, but
for couples:
$(t,m)$ is similar to $(t',m')$ if there exist  orientation-preserving
homeomorphisms of
$\varepsilon _i$ mapping each 
$I_c( t,m)$ to $I_c( t',m')$
  (the length of $\varepsilon _i$
  depends on
$m$, but this is irrelevant since we only consider the combinatorial
structure).  There are only finitely many similarity classes, and we
have seen in the previous proof that, as a non-metric tree, $T_t(m)$
depends only on the similarity class of $(t,m)$. This completes the
proof of Lemma \finit, hence also of the first assertion of Theorem
\contract{}.
\cqfd\enddemo

\subhead     Forester's deformation theorem\endsubhead

Before showing contractibility in the Gromov topology, we   use Lemma
\finitd{} to give a proof of Forester's deformation theorem [\Fo], which
we state as follows: 

\thm{Theorem \sta} Let $T_0,T$ be simplicial (non-metric) trees with 
the same elliptic subgroups. 
\roster
\item There exist trees    $T_0,T_1,\dots,T_k=T$ such that $T_{i+1}$
is obtained from
$T_i$ by an elementary expansion or an elementary collapse.  
\item If all edge stabilizers of $T_0$ and $T$ belong to a family $\A$ 
which is stable under taking subgroups, one may choose $T_i$ with edge
stabilizers in
$\A$. 
\endroster
\fthm

The second assertion is not in [\Fo], but   it can be
deduced from Forester's proof.

\demo{Proof} Using elementary collapses, we   assume that $T_0$ is 
reduced. Since elliptic subgroups of $T_0$ are elliptic in $T$, we can
map the vertex set of
$T_0$ equivariantly to the vertex set of $T$. Adjacent vertices have
distinct images by Remark \pasecr, and we may choose metrics on 
$T_0$ and $T$ so that there exists a morphism $f:T_0\to T$. 

   We apply Lemma \finitd{} (note that in its proof we only used the
fact that $T_0$ and $T$ have the same elliptic subgroups). The set of
intermediate trees
$T_t $ is contained in a  finite union of cones.
By Proposition \memetop, $T_t $ is a continuous path in the weak
topology. In particular, the union of all cones which meet this path
is connected in   that topology. Thus, there is a finite sequence of
closed cones
$T_0\in\C_0,\dots,\C_n\ni T$ such that for all $i\in\{1,\dots,n\}$, either
$\C_i$ is a face of $\C_{i-1}$ or  $\C_{i-1}$ is a face of $\C_i$.
Now if $\C'$ is a face of $\C$,   trees of $\C'$ are obtained from  
trees of $\C$ by an elementary collapse. This proves Assertion (1). 

Making $T_0$ reduced does not create new edge stabilizers. Edge
stabilizers of $T_t$ are in $\A$ because there is a morphism 
$\psi_t:T_t\to T$. This proves Assertion (2).
\cqfd\enddemo

\subhead Contractibility in the Gromov
topology\endsubhead

We now prove Assertion (2) of Theorem \contract. We assume that
$T_0$ has finitely generated vertex stabilizers, and we check that the
condition in Lemma \critere{} holds. We first wish to say that the
basepoint $P\in T$ depends continuously on $T$. As   $P$ does not
belong to a fixed space, we express it in the following way:

\nom\ptbas
\thm{Lemma \sta}
For any $a\in G$,  the map $T\mapsto d (P,aP)$, from $\D$ to $\R$,
is continuous. 
\fthm

Of course, $\D$ is now equipped with the Gromov topology.

\demo{Proof} We first consider the case when $\D$ is irreducible. We
recall the definition of $P$. We fix $g,h\in G$ with $g,h,[g,h]$ hyperbolic
in $\D$.  If the axes
$A_g$ and $A_h$ are disjoint,   $P$ is the point of $A_g$
closest to $A_h$. 
If $A_g$ and $A_h$ meet, $P$ is the largest element of the segment
$A_g\cap A_h$, with   $A_g$ ordered so that $g$ acts as a
positive translation. 

Using the general formula $d(x,gx)=\ell(g)+2d(x,A_g)$, it is easily
checked that
$P$ is characterized by the system
$$\left\{\aligned
d(P,gP)&= \ell(g) \\
d(P,hP)&= \ell(h)+2D \\
d(gP,hP)&=\ell(g)+\ell(h)+2D \\
d(gP,h\mi P)&=\ell(g)+\ell(h)+2D ,
\endaligned\right.\tag1$$
with $D$ equal to the distance between $A_g$ and $A_h$ (see Figure 2). 

\midinsert
\centerline{\includegraphics[scale=1]{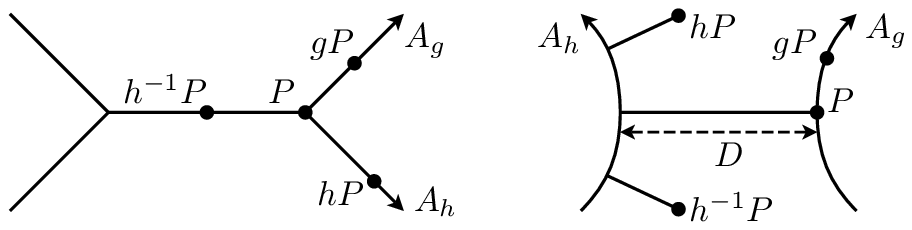}}
\botcaption{Figure 2}{The basepoint   $P$}
\endcaption  
\endinsert

As in the proof of Proposition \AbelianHausdorff, this system is stable:
if $\wtilde P$ satisfies each equation up to $\varepsilon $, then it is
$C\varepsilon $-close to $P$, for some fixed number $C$. Furthermore,
the right-hand sides depend continuously on $T$ (see [\Pa,
Proposition  3.5] for  the continuity of $D$).

We can now prove that $d (P,aP)$ is continuous at any given $T$.  Fix
$\varepsilon >0$.  Consider another $T'\in \D$, and the corresponding
system $(1')$.  
  If $T'$ is close enough to $T$ in the Gromov topology, then   there
exists a lift 
$\wtilde P \in T'$ satisfying $(1')$ up to $\varepsilon $. We
also require that $d(\wtilde P, a\wtilde P)$ be $\varepsilon $-close to
$d(P,aP)$. Then $\wtilde P$ is
$C\varepsilon $-close to the basepoint $P'$ of $T'$, so $d(P',aP')$ is
$(2C+1)\varepsilon
$-close to
$d(P,aP)$. 

When $\D$ is abelian, we  fix two 
positive hyperbolic elements $g,h\in G$ with distinct axes, and $P$ is
the endpoint of the ray
$A_g\cup A_h$. The argument is the same as in the irreducible case,
using the   system $(1_{gh})$ introduced in the proof of Proposition 
\AbelianHausdorff.
\cqfd
\enddemo

More generally:

\nom\ptbass
\thm{Lemma \sta}
Let $H_1,H_2$ be (possibly equal, possibly   trivial) finitely generated subgroups
of
$G$ which are elliptic in $\D$. Given $T\in\D$, let $P_i\in T$ be the projection of the
basepoint $P$ onto the fixed point  set of $H_i$. For any
$a\in G$,  the map
$T\mapsto d (P_1,aP_2)$  is continuous. 
\fthm

\demo{Proof}  We show continuity  at a given   $T\in \D$. For
$i\in\{1,2\}$, fix a finite generating set
$S_i$  for
$H_i$.   Choose $s_i
$  in $S_i$ so that $d(P,s_i P )$ is maximal  ($s_i$ is not necessarily 
unique). The point  $P_i$ is the midpoint of  
$[P,s_iP]$.

\midinsert
\centerline{\includegraphics[scale=1]{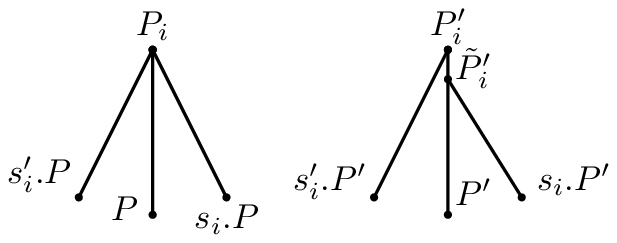}}
\botcaption{Figure 3}{Proof of Lemma \ptbass}
\endcaption  
\endinsert

Now consider   $T'$ close to $T$  in $\D$. Define  $P',P'_i,s' _i$
as in
$T$  (see Figure 3).  Note that $s_i$ is fixed,   but
$s'_i$ may vary as    $T'$ varies. Let
$\wtilde  P'_i$ be the midpoint of
$[P', s_iP']$.  Since $P'_i$ is fixed by $s_i$,  the point 
$\wtilde P'_i$    is located on $[P',P'_i]$, and
$d(\wtilde P'_i, P'_i)=\frac12(d(P', s'_iP')-d(P',s_iP'))$.  Our choice  of
$s_i$ guarantees  $d(P,sP)-d(P,s_iP)\le0$ for every $s\in S_i$.
Lemma
\ptbas{} therefore implies that $d(\wtilde P'_i, P'_i)$
 goes to
$0$ as
$T'\to T$.

Since distances between the four  points  $P ,
s_1P ,aP , as_2P $ vary continuously  by Lemma \ptbas, the  distance
$d(\wtilde P'_1,a\wtilde P'_2)$ between the midpoints of $[P',
s_1P']$ and $[aP', as_2P']$ converges to $d(P_1,aP_2)$ as $T'\to
T$. Thus $d( 
P'_1,a P'_2)\to d(P_1,aP_2)$.
\cqfd\enddemo

Since $f_T$ was defined using projections onto fixed point sets of
vertex stabilizers of $T_0$, this lemma immediately implies that the
condition in Lemma
\critere{} holds if $T_0$ has finitely generated vertex groups. This
completes the proof of Assertion (2) of Theorem \contract.

\subhead The closure of an irreducible space is
contractible\endsubhead

 We now wish to extend the contraction to
the closure $\ov\D$ of an irreducible
$\D$ (in the Gromov topology). Fix $g,h$  with $g,h,[g,h]$ hyperbolic.
See Proposition
\obord{} for properties of \Rt s $T\in\ov\D$. In particular, $T$ is 
irreducible. The characteristic sets $A_g$, $A_h$ (axis or fixed point
set) contain no tripod, and their intersection has length
$\le\ell(g)+\ell(h)$.

We have to change the definition of $P$, as there is no natural way to
orient $A_g$ if $g$ is elliptic. We use the ``symmetric'' definition: 
$P$ is the midpoint of the segment $A_g\cap A_h$ if $A_g$ and $A_h$
meet,   the midpoint of the ``bridge'' between them if they  are
disjoint.

\midinsert
\centerline{\includegraphics[scale=1]{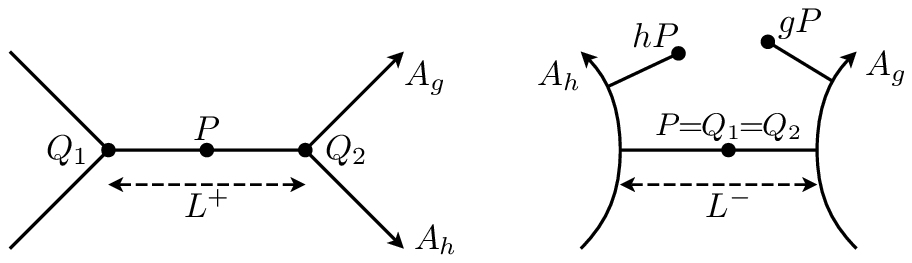}}
\botcaption{Figure 4}{The ``symmetric''     basepoint}
\endcaption  
\endinsert

The point  $P$ is characterized by the following
property (see Figure 4): there exist $Q_1,Q_2$  such that
$$\left\{\gathered d(P,gP)=d(Q_1,gQ_1)=d(Q_2,gQ_2)=\ell(g)+L^- \\
d(P,hP)=d(Q_1,hQ_1)=d(Q_2,hQ_2)=\ell(h)+L^- \\
d(P,Q_1)=d(P,Q_2)=L^+/2\\ d(Q_1,Q_2)=L^+ ,
\endgathered\right.\tag2$$
where $L^+=\max(L,0)$ and
$L^-=\max(-L,0)$, with 
$L\in\R$ equal to the length of $A_g\cap A_h$ if
$A_g\cap A_h$ is nonempty, to minus the distance between
$A_g$ and
$A_h$ if they are disjoint.

To see this, first suppose
$A_g\cap A_h=\ev$ (so   $L^+=0$). The equations reduce to
$d(P,gP)=\ell(g)+L^-$ and
$d(P,hP)=\ell(h)+L^-$ (with $Q_1=Q_2=P$). 
They   mean that
$P$ has distance
$L^-/2$ to both $A_g$ and $A_h$. This uniquely  determines it. 
If  $A_g\cap A_h\neq\ev$, the first two
equations mean that $P,Q_1,Q_2$ belong to $A_g\cap A_h$, the other 
two that the segment 
$[Q_1,Q_2]$ has length $L^+$ and midpoint $P$. This uniquely 
determines $P$.  

Furthermore, the system (2) is stable in the   sense
used above: 
 if 
  $(\wtilde P, \wtilde Q_1,\wtilde Q_2) $ is an approximate solution,
then   $\wtilde P$ is close to $P$.  

To extend Lemmas \ptbas{} and \ptbass{} to $\ov\D$, we have to know
that $L$ is continuous on $\ov\D$. It is continuous on $\D$ by [\Pa,
Proposition  3.5], but the argument of [\Pa]
breaks down when both
$g$ and $h$ are elliptic in $T\in\ov\D$. If this happens, first suppose that
the fixed point sets of $g$ and $h$ are disjoint. Then we have
$L=-\frac12(\ell(gh)-\ell(g)-\ell(h))$ in $T$ and in nearby trees, so $L$
is continuous at $T$. If $g$ and $h$ have a common fixed point, we know
that it is unique, so $L=0$ at $T$. Continuity at $T$ follows from the
general inequality
$|L|\le\max (\ell(g)+\ell(h),\frac12\ell(gh))$. 

The rest of the argument is now as before. The space  
$\ov\D$ is contracted   into a closed cone, as certain edges 
 of $T_0(T)$ may
have length $0$ if $T\in\ov\D\setminus\D$. 

\head   \sect Non-ascending deformation spaces  \endhead

Given a $G$-tree $T$, we now view the  quotient graph $\Gamma
=T/G$ as a graph of groups.  We often view a deformation
space as a collection of marked graphs of groups (a marking of
$\Gamma $ is an isomorphism from $G$ to the fundamental group of
$\Gamma $, well-defined up to composition with inner
automorphisms). We say that
$\Gamma
$ is reduced if
$T$ is. 

We use the same notation (such as   $v$)  to denote a vertex of $T$ or 
of $\Gamma $. The corresponding group (vertex stabilizer or vertex
group) is then denoted by
$G_v$.  Edges are denoted by letters such as $e$, with group $G_e$. An 
edge of
$\Gamma
$ is a {\it loop\/} if its endpoints are equal, a {\it segment\/} otherwise. 

\nom\slides
\subhead Slide moves\endsubhead

Let $e=vw$ and $f=vx$ be adjacent edges in a $G$-tree $T$, with
$G_f\inc G_e$. Assume   that $e$,  $f$ are not in the same orbit
as non-oriented edges. We define a new tree $T'$ with the same vertex
set as $T$ by replacing $f$ by an edge
$f'=wx$ and extending this operation equivariantly to the orbit of $f$.
We say that {\it $f$ slides across $e$\/} (see [\FoGBS]).  More generally,
one can slide several edges
$f_i$ satisfying $G_{f_i}\inc G_e$. Though sliding makes sense for
metric trees, we will apply it only in the context of non-metric trees.

The trees $T$ and $T'$ belong to the same deformation space (they
have the same vertex stabilizers). In fact, a slide move may be viewed
as an elementary expansion followed by an elementary collapse (in the
middle tree, $f$ is attached to the midpoint of
$e$). The trees $T$ and $T'$ even have the same edge stabilizers: since
$G_v\cap G_x\inc G_v\cap G_w$ and
$G_x\cap G_w\inc G_v$, we have $G_{f'}= G_w\cap G_x=G_v\cap
G_x=G_f$.

It is   easier to visualize a slide move on
$\Gamma $, as one edge sliding across another (one or both may be
loops). The number of vertices and edges of   $\Gamma $ does not
change. We usually describe a slide  move by its action on
$\Gamma
$. There may exist different moves in $T$ with the same projection in
$\Gamma $, but this will not be an issue here (compare [\Clb]).

\subhead Non-ascending spaces: definition and examples\endsubhead

Given an edge $e=vw$ of $\Gamma $ (possibly with $v=w$), there are
injections $G_e\hookrightarrow G_v$ and $G_e\hookrightarrow G_w$.
We attach two {\it labels\/} to the edge $e$, one near $v$ and one near
$w$. Each of these labels is $=$ or $\neq$, depending on whether the
corresponding injection is onto or not.  

With these notations, $\Gamma
$ is reduced if and only if all   labels $=$  are carried by loops.  If a
vertex $v$ has degree 2, at least one of the two labels near $v$ is 
$\ne$ (because there is no redundant   vertex in $T$). If
$e$ carries a label $=$ near $v$, then all  {other} edges attached to $v$
may slide across $e$.

When viewed in $\Gamma $, an elementary collapse (see Section   3)
is the collapse of a segment $e$ carrying at least one  $=$ label; such
an  
$e$ is called collapsible.  The graph $\Gamma $ is reduced if and only if
it contains no collapsible edge.   All collapses considered here will be
elementary collapses, so we often drop the word
elementary.
 
We usually identify the edges of the collapsed graph with the edges of
$\Gamma
$ other than $e$, and a vertex of $\Gamma $ with its image in the
collapsed graph. When
$e$ is collapsed, a collapsible edge $f\ne e$ may become
non-collapsible; a non-collapsible edge remains non-collapsible.

A {\it   strict ascending loop\/} [\FoGBS] is a loop $e$ of $\Gamma $
carrying exactly one $=$ and one
$\neq$.  The fundamental group of $e$, viewed as a graph of groups, is
then a strictly ascending HNN extension. All edges adjacent to the
basepoint of
$e$ may slide around $e$. 

A deformation space
$\D $ is {\it \na{}\/} if it is irreducible, and no  
$\Gamma
$ in
$\D$ contains a strict ascending loop. Note that  genuine abelian spaces
always contain graphs with a strict ascending loop.

The following observation will be very useful. Define an {\it ascending
circle\/} as an embedded  circle
$C\inc\Gamma $ which may be oriented in such a way that all edges of
$C$ carry a
$=$ label near their origin. If $\Gamma $ belongs to a \na{} space, then  
edges of $C$  also carry a $=$ near their terminal  point: otherwise,
collapsing all edges of $C$ except  one   would create a strict
ascending loop.

\thm{Proposition \sta} Let $\D$ be an irreducible deformation space. If
one of the following conditions is satisfied, then $\D$ is \na:
\roster
\item No   {\rm reduced} graph of groups in $\D$ contains a strict
ascending  loop. 
\item There exists a tree $T\in\D$ such that no fixed point set
$\fix(g)$, for
$g\neq1$, contains an infinite ray.
\item No group in $\A_{min}$ properly contains a conjugate of itself.
\item There exist $T_c\in\D$ such that, for any $T\in\D$, there exists
$\wtilde T\in\D$ such that both $T_c$ and $T$ may be obtained from
$\wtilde T$ by   collapses. 
\endroster
\fthm

Recall that the weak acylindricity condition (2) holds for all trees in
$\D$ if it holds for one (see Section 4), and that $\A_{min}$ consists of
groups fixing an edge in   every $T\in\D$ (see Definition \aminn).
Condition (3) holds if $\D$ is the canonical deformation space of a GBS
group which does not contain any $BS(1,n)$ with $n\ge2$ (see [\LeGBS]).
Condition (4) means that
$T_c$ is ``compatible'' with every
$T\in\D$: the graphs of groups $\Gamma _c$ and $\Gamma $ have a common
refinement $\wtilde\Gamma  $ in $\D$.

\demo{Proof} Starting with $\Gamma $ having a strict ascending loop
$e$ based at
$v$, we make it reduced without destroying the loop. First we collapse
all collapsible edges which don't contain $v$, or have a $\neq$ label
near $v$, or have two $=$ labels. This does not change $G_v$, so $e$
remains strictly ascending. If there remains  a   collapsible edge $f$,  
we change its $=$ label   near $v$ to
$\neq$ by sliding it around
$e$. This makes $f$ non-collapsible. Iterate until all  edges have been
made non-collapsible.

If $\Gamma $ has a strict ascending loop $e$, the associated   tree $T$
does not satisfy (2), and $G_e$ properly contains a conjugate 
($G_e\in\A_{min}$ if $\Gamma $ is reduced). 

There remains to  consider (4). Recall that   length functions of trees
in
$\D$ are bi-Lipschitz equivalent (Theorem \lessix). To estimate the
bi-Lipschitz constant, we restrict to trees all of whose edges have
length $1$. Suppose $T'$ is obtained from $T$  by collapsing all edges in
the orbit of an edge 
$e$. As pointed out in [\Fo], components of this orbit have diameter at
most 2.  This implies $\ell_T/3\le \ell_{T'}\le \ell_T$ (at most two
thirds of the edges in a given axis may be collapsed).
  If we can pass from $T$ to $T'$ by $p$ elementary operations
(expansions or collapses), we get
$3^{-p}\ell_T\le\ell_{T'}\le3^{p}\ell_T$. 

In the situation of (4), the number of operations needed to pass from
$T_c$ to $T$ is bounded by   $ |\Gamma _c|+|\Gamma  | $, where
$|\Gamma |$ denotes the number of edges of the quotient graph, as we
may assume that no edge of $\wtilde T$ gets collapsed in both $T_c$
and
$T$. 

A consequence of (4) is therefore the following: {\it Given $k$,
there exists $C$ such that, if $T,T'\in
\D$ have edges of length 1, and the quotient graphs $\Gamma ,\Gamma
'$ have at most 
$k$ edges, then $\frac1C\ell_T\le\ell_{T'}\le C\ell_T$.}
We now suppose that some reduced $\Gamma \in\D$  has a strict
ascending loop $e$ based at a vertex $v$, and   we obtain a
contradiction to this assertion. 

Let
$f$ be another edge adjacent to $v$ (it exists because $\D$ is
irreducible).  We may assume that it carries no $=$ label (if it does,
then it is a loop at $v$ because
$\Gamma $ is reduced, and we get rid of
$=$ labels by sliding around $e$). 

Define a graph $\Gamma _n$, with Bass-Serre tree $T_n$, by sliding $f$
$n$ times around $e$ (oriented from its $=$ label to its $\ne$ label).
Recall that all   trees $T_n$ have the same vertex set. Consider a lift
$v_0x$ of
$f$ to
$T_0=T$, with
$v_0$ lifting
$v$. It is slid to a position $v_nx$ in $T_n$.  For $0\le i\le n$, the
distance between
$v_n$ and
$x$ in
$T_i$ is $n-i+1$ (the geodesic is $v_nv_{n-1}\dots v_ix$). 

Choose $g\in G$ fixing $x$ but not $v_0$ (it exists because $f$ carries
no $=$), and
$h_n$ fixing
$v_n$ but not $v_{n-1}$ (the sequence $G_{v_n}$ is strictly ascending
because $e$ is). The distance between the fixed point sets of $g$ and
$h_n$ in $T_i$ is
$d_{T_i}(v_n,x)=n-i+1$, so
$\ell_{T_i}(gh_n)=2(n-i+1)$. Writing this for $i=0$ and $i=n$, we see 
that the bi-Lipschitz constant between $\ell_T$ and
$\ell_{T_n}$ goes to infinity with $n$,  contradicting the
consequence of (4) stated  above. 
\cqfd\enddemo

\subhead  Connecting trees by slide moves
\endsubhead

In the rest of this section, we only consider   \na{} spaces.

It is easy to show that, in a 
 \na{} deformation  space, every elliptic subgroup is contained in a
maximal one.   If
$T$ is reduced, there is a bijection  between the set of   vertices  of
$\Gamma $ and the set  of conjugacy classes of maximal elliptic
subgroups (this is the vertical set $\M$ of Section 4).  The  next result
will provide more precise information.  We consider trees as
simplicial, non-metric, objects (or we restrict to trees with edges of
length 1). 

\nom\slides
\thm{Theorem \sta} Let $\D$ be a \na{} deformation space. Any two
reduced simplicial trees
$T,T'\in\D$ may be connected by a finite sequence of  slides. 
\fthm

Conversely,  performing a slide move on a reduced tree in a \na{} space
  yields a reduced tree.

This theorem was proved by Forester [\FoGBS] for the canonical
deformation space of GBS groups. His proof works in the general
situation. We give a different argument. 

\demo{Proof}  By  Remark \passage, we may assume that $T,T'$ are
different reductions of a tree $S$. We also assume that no edge of $S$
gets collapsed in both
$T$ and $T'$. We work with the quotient graphs of groups $\Gamma
,\Gamma ',\Theta
$
  of $T,T',S$. We first show that there exist edges
$e,e'$ of $\Theta $ with a common vertex $v$, such that both $e$ and $e'$
carry a  $= $ near $v$, and $e$ gets collapsed in $\Gamma '$ and $e'$
gets collapsed in $\Gamma
$.

Since $\Theta $ is not reduced, we can find an edge $e_1$ of $\Theta $
which gets collapsed in $\Gamma '$. Let $v_1$ be an endpoint of $e_1$
with label $=$.  Since $e_1$ is not collapsed in $\Gamma
$,    some edge $e'_1$ of $\Theta $ with   endpoints $v_1, v'_1$ gets
collapsed in
$\Gamma
$.  If the label of $e'_1$ at $v'_1$  is $\ne$, its label at $v_1$ is $=$ and
we are done. If not, we repeat the argument, using $e'_1,v'_1$ instead of
$e_1,v_1$. We obtain edges $e_2,e'_2,e_3,\dots$ which get collapsed
in  
$\Gamma $ or $\Gamma '$ alternatively. If the process does not stop, we
eventually create an ascending circle $C$. Since $\D$ is \na, $C$ carries
only $=$ labels and we can find $e,e'$.

\midinsert
\centerline{\includegraphics[scale=1]{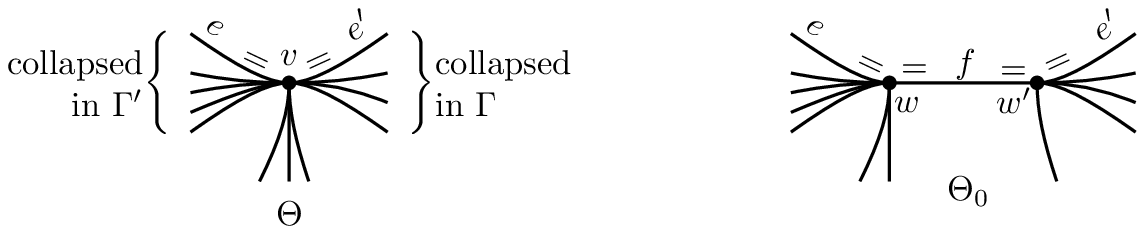}}
\botcaption{Figure 5}{Proof of Theorem \slides{}}
\endcaption   
\endinsert 

Let $e,e'$ be as constructed. We define a new graph of groups $\Theta
_0$ by performing an elementary expansion at $v$  in $\Theta
$, as follows (see Figure 5). We replace
$v$ by a new  edge $f=ww'$   carrying two $=$ labels.  Edges of $\Theta $
incident to
$v$ are attached   to   $w$ or   $w'$, the only restriction being that   
edges collapsed in
$\Gamma '$ (resp\. $\Gamma  $) are   attached to   $w $ (resp\. $w '$)
(changing the attachment point of an edge  amounts to sliding it across
$f$).  In particular, the edge
$e$ (resp\. $e'$) is attached to $w$ (resp\. $w'$). It carries a label $=$ at
$w$ (resp\.
$w'$). 

 Let $F_0$ be the image in $\Theta _0$ of the subforest of
$\Theta
$ which gets collapsed in $\Gamma  $ (it contains
$e'$). Collapsing $F_0\cup
f$ in $\Theta _0$ yields $\Gamma $. 
Since no edge of $F_0$ is attached to  $w $, we may also collapse 
$F_0\cup e$. This gives a graph $\Gamma _1$ which differs from $
\Gamma $ by slide moves (the middle graph is obtained from $\Theta
_0$ by collapsing $F_0$). Define
$\Gamma '_1$ similarly, using $\Gamma '$ instead of $\Gamma $. The 
graphs
$\Gamma _1$ and $\Gamma '_1$ are both reductions of the graph 
$\Theta _1$ obtained from
$\Theta _0$ by collapsing $e\cup e'$. Since $\Theta _1$ has one less 
edge than
$\Theta $,    iterating this construction
yields a sequence of slides connecting $\Gamma $ and $\Gamma '$.
\cqfd\enddemo 

\nom\memgrp
\thm{Corollary \sta} If $\D$ is \na, all reduced trees in $\D$ have the
same vertex and edge stabilizers; given two reduced trees in $\D$,
there are   $G$-equivariant bijections
$\psi _V$, $\psi _E$ between their vertex  and   edge sets. \cqfd
\fthm

By  [\Di], the existence of  $\psi _V$ always implies that of $\psi _E$. In
general,
$\psi _V$ and
$\psi _E$ are not compatible with the incidence relations.

\subhead A finite-dimensional retract \endsubhead

Recall that we only consider  graphs of groups in \na {} deformation
spaces. All collapses considered here are  elementary.

As mentioned above, a collapsible edge of $\Gamma $ may become
non-collapsible when other edges are collapsed.  This motivates the
following definition. 

An edge $e$ of $\Gamma $ is a {\it surviving edge\/} if it is
non-collapsible, or may be made non-collapsible by collapsing other
edges. In other words, $e$ is a surviving edge if and only if one can
make $\Gamma $ reduced without collapsing
$e$.

  Recall that there are two types of non-collapsible edges :  
segments carrying two
$\neq$ labels, and loops.  In order  to describe surviving
edges,  we define  two types of
subgraphs
$\gamma $ of $\Gamma$, in terms of their topology and the labels they
carry:

$\bullet$   $\gamma =(e_1,\dots, e_k)$ is an embedded edge-path
(homeomorphic to
$[0,1]$); all labels carried by $\gamma $ at interior vertices are $=$,
but the labels carried by $e_1$ at its origin, 
    and  by $e_k$ at its terminal  point,  are $\neq$ (when $k=1$, this
is just a non-collapsible segment). 

$\bullet$  $\gamma $ is homeomorphic to a circle;   all  labels carried by
$\gamma
$ are $=$, except possibly those  at a  special vertex $v_0$  (if there is
a special vertex,   both labels at $v_0$ are $\neq$ by the
\na {} assumption).

We call such a $\gamma $ a {\it shelter\/}. Note that one-edge shelters
are exactly non-collapsible edges. Any edge
$e$ contained in a shelter $\gamma $ is a surviving edge, since it
becomes non-collapsible when the other edges of
$\gamma
$ are collapsed. The following lemma will imply the converse (every
surviving edge is contained in  a shelter).

\nom\shelt
\thm{Lemma \sta} Let $\Gamma $ be a graph of groups belonging to a
\na{} deformation space, and let
$\Gamma '$ be obtained from $\Gamma $ by collapsing a collapsible
edge
$f$. An edge $e\neq f$ is contained in a shelter of $\Gamma $  if and only
if its image
$e'$ in $\Gamma '$ is contained in a shelter.
\fthm

\demo{Proof  } This is just a statement about graphs with a $=,\neq$
labeling.  The proof is not difficult, but there will be several cases to
consider. 
 
Let $v,w$ be the vertices of
$f$. They are distinct and
$f$ carries at least one $=$ label, say near $v$. When $f$ is collapsed,
labels carried by other edges near  
$w$ don't change (the vertex group does not change). A $\neq$ label
near $v$ remains $\neq$.  A $=$ label gets replaced by the label of $f$
near $w$. 

Assume that $e$ is contained in a shelter $\gamma $.   If  $f$ is  disjoint
from $\gamma $, or contained in $\gamma $, or if $f\cap\gamma
=\{w\}$, the image of
$\gamma
$ in
$\Gamma '$ is a shelter. If $f\cap\gamma =\{v\}$,  the image of
$\gamma  $ consists of one or two shelters. Now assume
$f\cap\gamma =\{v,w\}$. If the label of $f$ near $w$ is $\ne$, the \na{}
hypothesis implies that 
$w$  must be an endpoint of $\gamma $, or the special vertex $v_0$.
The image of
$\gamma $ in $\Gamma '$ is homeomorphic to a circle, or a wedge of
two circles, or the union of a circle and an interval. One checks that
each of these circles or intervals is a shelter.  We have proved one
direction of the lemma.

Conversely,  assume that $e'$ is contained in a shelter 
$\gamma '$. Let $w'$ be the image of $v$ and $w$ in $\Gamma '$  (so  
$\Gamma $  is obtained from
$\Gamma '$ by expanding  
$w'$ into $f$). We may assume $w'\in\gamma '$. 
 Let $\gamma _0$ be the subgraph of  
$\Gamma $ consisting of edges projecting onto the edges of $\gamma '$
(it has one or two components). Edges of
$ \gamma _0$ near $f$ may be attached to either $v$ or $w$. We may
assume that at least one is attached to $v$ (otherwise, $\gamma _0$ is
obviously a shelter).

 First assume that  the label of $f$ near $w$ is $=$. The only nontrivial
case is when $\gamma _0$ has one edge attached to $v$ and one to $w$. 
If $w'$ is the special vertex  of a circular $\gamma '$, then $\gamma _0$
is a shelter (homeomorphic to an interval).  If not, then
$\gamma _0\cup f$ is a shelter.

Now suppose that the  label of $f$ near $w$ is $\neq$.
 Then $w'$ has to be an endpoint of $\gamma '$, or its special vertex
$v_0$, since otherwise no edge of $\gamma _0$ could be attached to
$v$. If it is an endpoint, then $\gamma _0$ is a shelter if the label of
$\gamma _0$ near $v$ is $\ne$, and   $\gamma _0\cup f$ is a shelter if
the label is
$=$. The same conclusion holds if $w'$   is the special vertex of
$\gamma '$ and
$\gamma _0$ has an edge attached to $w$. If  $w'=v_0$, and
$\gamma _0$ has two edges attached to $v$, then $\gamma _0$ is a
shelter.
\cqfd\enddemo

\nom\surv
\thm{ Corollary \sta} Let $\Gamma $ be a graph of groups belonging
to a \na{} deformation space. \roster
\item An edge is a surviving edge if and only if it is contained in a
shelter.
\item  Let
$\Gamma '$ be obtained from $\Gamma $ by collapsing a collapsible
edge
$f$. An edge $e\neq f$ is a surviving edge of $\Gamma $  if and only if
its image $e'$ in $\Gamma '$ is a surviving edge.
\endroster
\fthm

\demo{Proof} We have seen that edges contained in a shelter are
surviving edges. If $e$ is a surviving edge, one may collapse edges
other than $e$ to obtain a reduced $\Gamma _r$. The image of $e$ in
$\Gamma _r$ is non-collapsible, so is a one-edge shelter. By the ``if''
part of the lemma,  
$e
$ is contained in a shelter of $\Gamma $. 
  Assertion (2) is just a rewording of the lemma.
\cqfd\enddemo

Let $\G\inc\D$ be
the set of graphs of groups all of whose edges are surviving edges. It is
stable under collapses, so its projection $P\G\inc P\D$ is a subcomplex.
Using shelters, it is easy to see that trees corresponding to graphs
$\Gamma
$ in
$\G$ are BF-reduced (as defined in Section 3):   if a vertex $v$ of
$\Gamma $ has valence 2 and the labels near $v$ are $=,\neq$, the edge
carrying the $=$ label cannot be contained in  a shelter.
 Proposition
\findim{} implies that
$\G$ is finite-dimensional.

\nom\spine
\thm{Theorem \sta} Let $\D$ be a \na{} deformation space, equipped
with the weak topology. There is a
  deformation retraction $r$ from $P\D$ onto the   
finite-dimensional subcomplex $P\G$. If  
$F\inc\Out(G)$ leaves $\D$ invariant, then
$r$ is $F$-equivariant.
 \fthm

\demo{Proof} 
The deformation retraction  is simply given
by letting the length of the non-surviving edges go linearly to $0$. It
maps a   given closed simplex into itself in a continuous way by
Assertion (2) of Corollary \surv,
 so it is continuous in the weak topology. The deformation is clearly
equivariant with respect to the group of automorphisms leaving $\D$
invariant (its action on a graph  of groups consists in changing the
marking).
\cqfd\enddemo

 We   have seen that all reduced trees in $\D$
have the same edge stabilizers (Corollary \memgrp). The same 
property holds for trees in $\G$, since all edges of graphs in $\G$ are
surviving edges. This shows that $\G$ is contained in the reduced
deformation space $\D_r$ defined in Section 4. 

\example{Examples }  In these examples, $G$ is a GBS group, $\D$ is
its canonical deformation space (see  the end of Section 3), and
$\G$ is the set of graphs of groups all of whose edges are surviving
edges. 

$\bullet$ Let
 $G =\langle
a,b,c,t\mid a^2=b^6, b^2=c^2, tc^3t\mi=b^9\rangle$. Using arguments  
from [\FoGBS], it may be shown
  that
$\D$ is
\na{} (although $G$ contains $BS(1,3)$). The subcomplex
$\G$ is 4-dimensional, but 
$\D_r$ is infinite-dimensional.

$\bullet$ Let $G=BS(2,4)$. In this case, $\G$ is not a subcomplex and
there is no retraction from $\D$ to $\G$.  Clay has shown, however, that
there always is a retraction from $\D$ to the barycentric spine of $\G$
[\Clb]. 

$\bullet$ Let
 $G =\langle
a,s,t\mid  sa s\mi=a^2,  ta t\mi=a^2\rangle$. By sliding one loop
many times around the other and performing elementary expansions,
one constructs  graphs of groups with an ascending circle of arbitrarily
large length. It follows that $\G$ is infinite-dimensional.

\endexample

\head \sect Applications to automorphisms \endhead

\subhead  Restriction   \endsubhead

Up to now, we have been working with a single finitely generated
group
$G$. In this subsection, we also  consider a  finitely generated group
$\widehat G$ containing   $G$
 as a normal subgroup.   
This study will be   used in the next subsection, with 
  $G$ a centerless group and $\widehat G$   the preimage
of some $F\inc\Out(G)$ under the natural map $\Aut(G)\to\Out(G)$.

The results of the present subsection also hold (and
are easier to prove)   if
$G$ has finite index in
$\widehat G$ (without being normal) and more generally if
 any $h\in \widehat G$ has a power in $G$.

We will consider irreducible trees.  Recall that they are determined by
their length function $\ell$.  The Gromov topology and  the axes
topology coincide on the set of irreducible trees (see Section 5). In
this   section it will be more convenient to use the axes topology
(i.e\. to work with length functions). 

  Any $\widehat G$-tree $T $ may be viewed as a $G$-tree. We call this
a {\it restriction\/}.   We sometimes  write $T_{\widehat G}$ and
$T_G$ instead of simply $T$ to emphasize which group we're
considering. Note that
$T_{\widehat G}$ is obviously minimal and irreducible if $T_G$ is. 

\nom\formu
\thm{Lemma \sta} Let $G\normal \widehat G$ be a pair of finitely 
generated groups with $G$ normal in $\widehat G$.  Let $T$ be an
irreducible
$\widehat G$-tree, with length function $\ell$. Assume that $G$ is  not 
contained in the kernel of the action.
\roster
\item The   $G$-tree  $T_G$ is minimal and irreducible.
\item For any $h\in \widehat G$, one has   $$\ell(h)=\sup_{g\in G}\    
\limsup_{k\to+\infty}\frac1{2k}\,\ell([g,h^k]).$$
\item $\ell $ is determined by its restriction to $G$.
\endroster
\fthm

 The kernel   of the action is the set of
elements acting on $T$ as the identity.

\demo{Proof}   If
$G$ is elliptic, its fixed subtree is
$\widehat G$-invariant because $G$ is normal, so equals $T$ by
minimality of
$T_{\widehat G}$; this means that $G$ is contained in the kernel. If
there is a unique $G$-invariant end, or a unique
$G$-invariant line, this end or line is   $\widehat G$-invariant,
contradicting irreducibility of $T_{\widehat G}$. The only remaining
possibility is that $T_G$ is irreducible. It is minimal, because the
unique minimal $G$-invariant subtree is $\widehat G$-invariant and
$T_{\widehat G}$ is minimal. 

Given any $h$, denote by $A_h$ its axis or fixed point set. 
Recall the formula $d(x,hx)=\ell(h)+2d(x,A_h)$, valid for any $x\in
T$, and the general inequality
$\ell(hh')\le\ell(h)+\ell(h') +2d(A_h,A_{h'})$, obtained by computing
$d(hh'x,x)$ for an $x\in A_{h'}$ such that
$h'x$ minimizes the distance to $A_h$. Writing $[g,h^k]=g(h^kg\mi
h^{-k})$, we obtain $$\ell([g,h^k])\le 2\ell(g)+2d(A_g, h^kA_g)\le
2\ell(g)+2k\ell(h)+4d(A_g,A_h).  $$
This shows that $ \ell(h)$ is bigger than the sup appearing in the
statement of the lemma.

To prove the opposite inequality, we may assume that $h$ is
hyperbolic. By Assertion (1), there exists a hyperbolic $g\in G$ with
$A_g\cap  A_h$ compact. For $k$ large, $A_g$ and $A_{h^kg\mi
h^{-k}}=h^kA_g$ are disjoint, and their distance is
$k\ell(h) +o(k)$. Since $\ell(gg')=\ell(g)+\ell(g') +2d(A_g,A_{g'})$ for
hyperbolic elements with disjoint axes, we have $\ell([g,h^k])=2k\ell(h)
+o(k)$.

Assertion (3) follows from (2),   since
$[g,h^k]\in G$ if $g\in G$.
\cqfd\enddemo

Let $G,\widehat G$ be as in the lemma. Consider any 
$\widehat G$-tree $T_{\widehat G}$, and its
restriction $T_G$. Elliptic subgroups of $T_G$ are precisely
intersections  of elliptic
 subgroups of $T_{\widehat G}$ with $G$. Restriction therefore sends a 
given deformation space
$\D_{\widehat G}$ of $\widehat G$-trees into a deformation space 
$\D_G$ of $G$-trees.  We denote this map by $i:\D_{\widehat G}\ra
\D_G$. By Lemma \formu, it is injective if $\D_{\widehat G}$ is
irreducible and $G$ is not elliptic in $\D_{\widehat G}$. 

\nom\gh
\thm{Lemma \sta} Let $G\normal \widehat G$ be a pair of finitely 
generated groups with $G$ normal in $\widehat G$.
\roster
\item An irreducible deformation space of $G$-trees contains  the
image of  at most one deformation space of $\widehat G$-trees. 
\item If $\D_{\widehat G}$ is an irreducible deformation space of 
$\widehat G$-trees in which $G$ is not elliptic,  restriction induces a
homeomorphism from
$\D_{\widehat G}$ onto a closed subspace of an irreducible 
deformation space
$\D_G$ of $G$-trees. 
\endroster
\fthm

Assertion (2) is true both in the  
Gromov (or axes) topology and in the weak topology.

\demo{Proof}  Recall that two
irreducible trees belong to the same deformation space if and only if 
their length functions are bi-Lipschitz equivalent (Theorem \lessix).  If
the restrictions of $\ell$ and $\ell'$ to $G$ are   bi-Lipschitz equivalent,
then $\ell$ and $\ell'$   are   bi-Lipschitz equivalent by Assertion
(2) of Lemma
\formu. This shows the first assertion.

Under the hypotheses of (2), we   have seen that restriction induces an
injection
$i:\D_{\widehat G}\to\D_G$ into an irreducible space. We now study the
topological properties of $i$.

First   consider the weak topology.
Given a
$\widehat G$-tree
$T_{\widehat G}$, restriction maps the closed cone spanned by
$T_{\widehat G}$ linearly  into the closed cone of $\D_G$ spanned by $T_{ 
G}$ (as in Section 5, we parametrize cones by edge-lengths). The image
consists of all trees satisfying certain equalities between edge-lengths.
In particular, at most one open cone of
$D_{\widehat G}$ maps into a given open cone of $\D_G$. Assertion (2) for
the weak topology follows from these observations. 

In the axes topology, the only obvious fact is continuity of $i$. To prove
that $i$ is a homeomorphism onto a closed set, we suppose that $\ell_n$ is
a sequence of length functions on $\widehat G$, associated to trees in
$\D_{\widehat G}$, whose restrictions to
$G$ converge to a nontrivial length function
$\ell_G$ associated to a tree $T_G$ in $\D_G$. 
We show that $\ell_n$ converges to a length function  associated to
a tree in $\D_{\widehat G}$.

Since the set of
projectivized length functions on
$\widehat G$ is compact [\CM], we can find numbers $a_n>0$ such that $a_n\ell_n$
converges to a nontrivial length function $\ell $ on $\widehat G$ (after replacing
$\ell_n$ by a subsequence). The sequence $a_n$ is bounded, because
$\ell_G$ is nontrivial. We have to bound it away from $0$.

 Let $T$ be the \Rt{} with length function $\ell$.   It 
is irreducible   by Proposition \obord.  If
$g\in G$ is hyperbolic in
$\D_{\widehat G}$, it may be elliptic in $T$, but cannot fix a tripod (see Proposition
\obord). This implies that $G$ is not contained in the kernel of the
action. In particular (see Lemma \formu), there exists $g_0\in G$ with
$\ell(g_0)\neq0$. Since
$\ell_n(g_0)$ converges to $\ell_G(g_0)$, and $a_n\ell_n(g_0)$
converges to a nonzero number, the bounded sequence $a_n$ has a
positive limit, which we may assume to be $1$. This implies that
$\ell_n$ converges to $\ell $, and 
$\ell_G$ is the restriction of $\ell$ to $G$. In particular, $T_G$ is the
restriction of $T$, and $T$ is simplicial. It is in $\D_{\widehat G}$
by Assertion (1). This completes the proof of  Assertion (2) for the axes  
topology.
\cqfd
\enddemo

\subhead  Fixed point sets \endsubhead

We return to the study of a single  finitely generated group $G$. 
There is a natural action of $\Out(G)$ on the set of $G$-trees, given by
precomposition.  If a subgroup of $\Out(G)$ fixes a tree, it acts on
the deformation space containing it.

\nom\fixcontr
\thm{Theorem \sta} Let $\D $ be an irreducible deformation
space.
If a  finitely generated group $F\inc\Out(G)$ fixes a tree $T_0\in\D $,
the fixed point set of
$F$ in
$\D $ is contractible in the weak topology. 
 \fthm

 \demo{Proof}  Let $\widehat G$ be the preimage of
$F$ by the natural map $\Aut(G)\to\Out(G)$. First assume that $G$ has
trivial center. It embeds into $\Aut(G)$, as inner automorphisms, and 
$\widehat G$ is a finitely generated group   containing
$G$   as a normal subgroup. 

The key observation is that a
$G$-tree $T_G$ is $F$-invariant if and only if the action of $G$ on
$T$ extends to an action of $\widehat G$. 
 The fixed point set of $F$ in $\D$ 
is homeomorphic
to an irreducible deformation space $\D_{\widehat G}$ by Lemma \gh,  
so  is contractible by Theorem
\contract. 

If $G$ has a nontrivial center, $\widehat G$ contains $G/Z(G)$ as a 
normal subgroup and we simply view 
$\D$   as a deformation space of $G/Z(G)$-trees (see the last example 
in Section 3).
\cqfd\enddemo

Theorem
\fixcontr{} also holds   in the Gromov topology, under the extra
assumption that $T_0$ has finitely generated vertex stabilizers (as a
$G$-tree, hence also as a $\widehat G$-tree).

\thm{Corollary \sta}  Let $\D $ be an irreducible deformation
space. Suppose that $\D$ is finite-dimensional, or \na.
If a solvable finite
group
$F\inc\Out(G)$ leaves $\D$ invariant, then $F$ has a fixed point in $\D$.
\fthm

\demo{Proof} By induction on the order of $F$. If $F$ is not trivial, it has
a normal subgroup
$F_1$ with
$F/F_1$ cyclic of prime order $p$. By induction, $F_1$ fixes a point in
$\D$. Let $\D_1\inc \D$ be its fixed point set.

First assume that $\D$ is finite-dimensional. Then 
 $\D_1$ is  
finite-dimensional, and contractible   by Theorem
\fixcontr. The group  
$F/F_1\simeq\Z/p\Z$ acts on $\D_1$, and   has a fixed point  
(otherwise, it would have a finite-dimensional classifying space). This
means that
$F$ has a fixed point in
$\D$. 
 
If $\D$ is \na, we consider the finite-dimensional subcomplex 
$\G\inc\D$ constructed in Section 7. The deformation retraction
$\D\to
\G$ is
$F$-equivariant.
Let $\D_1$ (resp\. $\G_1$)  be the  
fixed point set of $F_1$ in  $\D$ (resp\. $\G$). By induction, $\D_1$ is
nonempty. As it is  contractible, we deduce that   $\G_1$ is
contractible and we conclude by considering the action of $F/F_1$ on
$\G_1$. 
\cqfd\enddemo

 \subhead Locally finite trees\endsubhead

In this subsection, we consider a deformation space $\D$ consisting of locally finite
trees, and the projectivized space $P\D$. We let $\Out_\D(G)$ be the subgroup of
$\Out(G)$ leaving $\D$ invariant.

We say that $H\inc G$ is a stabilizer if it is a vertex or edge
stabilizer in some
$T\in\D$. All  stabilizers   are commensurable as
subgroups of $G$: if $H_1$,  $H_2$ are stabilizers, 
then $H_1\cap H_2$ has finite index in both $H_1$ and $H_2$ (this is 
clear if
$H_1,H_2$ are stabilizers in the same tree, and also if they are
stabilizers in trees differing by an elementary expansion).

A group $H$ (possibly infinitely generated) is small [\BFa] if there is no
$H$-tree in which axes of two hyperbolic elements intersect in a
compact set. Any group not containing
$F_2$ is small.  Smallness  is
not stable under taking subgroups, even of finite index.

To remedy   this, we say that $H\inc G$ is
\emph{small in $G$} if there is no $G$-tree in which axes of two
hyperbolic elements of $H$ intersect in a compact set. Of course, any
small group contained in $G$ is small in $G$, and being small in $G$
is stable under taking subgroups. It is also a commensurability
invariant.

\nom\alge
\thm{Lemma \sta} Given $G$, there exists at most one irreducible
deformation space consisting of locally finite trees   whose
stabilizers are small in $G$.
\fthm

\demo{Proof} This is proved in [\Fo]   in the case of stabilizers
$\Z$, in [\Clfp] for slender stabilizers (recall that a group is slender
[\DS] if all of its subgroups are finitely generated). We use the same
argument to determine elliptic subgroups algebraically: if
$T$ is locally finite, and
  stabilizers are small in $G$, then  a subgroup $H\inc G$ is elliptic if
and only if it is contained in a   subgroup which is small in $G$ and
commensurable to all its conjugates. 

The ``only if'' direction is obvious. 
Conversely, suppose $H$ is small in $G$, commensurable to all its
conjugates, but not elliptic. It is easy to check  that
$H$ acts on $T$ with  a unique fixed end or a unique invariant line. 
 This line or end is $G$-invariant, contradicting irreducibility. 
\cqfd\enddemo

A basic invariant of $\D$ is the {\it modular homomorphism\/}  $\Delta
:G\to \Q^{+*}$ with values in the multiplicative group of positive
rationals [\BK, \FoGBS]. It may be defined by choosing a stabilizer $H$
and setting $\Delta (g)=\frac{[H:H\cap gHg\mi]}{[gHg\mi:H\cap
gHg\mi]}$. Given $T\in\D$, one has $\Delta =\tau \circ \theta $,
where $\theta $ is the natural epimorphism from  $G$ to the topological
fundamental group of the quotient graph $\Gamma $ (see the  
beginning of Section 4), and
$\tau
$ is defined on a loop
$\gamma =(e_1,\dots,e_k)$ by $\tau (\gamma )=\prod_{j=1}^k\frac{
i(e_j)}{i(\bar e_j)}$, denoting by $i(e_j)$ (resp\. $i(\bar e_j)$) the index 
of the edge group   in the vertex group at the origin (resp\. endpoint)
of
$e_j$. 

\example{Remark}
There is a more refined invariant,  with values in the
abstract commensurator of $H$. For GBS groups, it is
the (signed) modular homomorphism used in [\FoGBS] or [\LeGBS].
\endexample

As in [\FoGBS], we say that $\D$ has  {\it no nontrivial integer
modulus\/} if the image of
$\Delta
$ contains no integer $n>1$ (when vertex groups are isomorphic to
$\Z$, this is equivalent to saying that $G$ does not contain a solvable
Baumslag-Solitar group $BS(1,n)$ with
$n>1$, see [\LeGBS]). An irreducible deformation space with no
nontrivial integer modulus is obviously \na.

\thm{Proposition \sta} Let $H$ be a finitely generated   subgroup of
$G$  such that $H$, and every group commensurable to $H$, has finite 
outer automorphism group. Let $\D$ be an irreducible deformation
space consisting of locally finite  $G$-trees with stabilizers
commensurable to $H$.  
\roster
\item If $H$ is small   in $G$, then $\Out_\D(G)=\Out(G)$.
\item If $H$, and all its finite index subgroups, have finite center, then
$\Out_\D(G)$  acts on $ \D$ with finite stabilizers.
\item If there is no nontrivial integer modulus, then  $\Out_\D(G)$ acts
on $P\D$ with finitely many orbits of simplices. 
\endroster
\fthm

\demo{Proof} (1) follows from Lemma \alge. 

For (2),  consider $T\in\D$. We study its stabilizer $\Out^T(G)\inc \Out
(G)$ using results and notations from   [\LeGD]. Since all edge and
vertex stabilizers have Out finite, the group of twists
$\T$ has finite index in $\Out^T(G)$ [\LeGD, Theorem 1.6]. This group
$\T$  is a quotient of a finite direct product of centralizers
$Z_{G_v}(G_e)$, with $G_e$ an edge stabilizer and
$G_v$ an adjacent vertex stabilizer. The hypothesis of (2)   implies
that all groups commensurable with $H$ have finite center, so all these
centralizers are finite.

  (3) was proved by Forester for GBS groups.  We generalize his
argument. Given $T\in\D$, consider the quotient graph of groups
$\Gamma $.  The hypothesis about moduli implies that there is a 
uniform bound $N$ (independent of $T$) for the index of an edge group
in an adjacent vertex group [\FoGBS, Theorem 8.2].

If $H_1,H_2$ are commensurable to $H$,
we note that  there are only finitely many injections 
$H_1\hookrightarrow H_2$ whose image has index $\le N$, up to
postcomposition with an inner automorphism of
$H_2$.  This is clear since $H_2$ has finitely many subgroups of index 
$\le N$, and
$\Out(H_1)$ is finite.

Since $P\D$ is  a locally finite complex  (see Section 5),
 and every closed simplex contains a
reduced tree, it suffices to show that there are finitely many
$\Out_\D(G)$-orbits of non-metric reduced trees.  
As $\D$ is \na, all reduced trees in
$\D$ are related by slide moves (Theorem \slides), so there are 
finitely many possibilities for $\Gamma =T/G$ as a graph. 
We now have to consider
$\Gamma  $ as a graph of groups, with edge groups, vertex groups, and 
inclusions.

All reduced   trees
  have the same edge and vertex stabilizers (Corollary \memgrp), so
there are finitely many possibilities for  the isomorphism type of   
edge and vertex groups of   $\Gamma $. As for the inclusions, we have
observed that there are finitely many possibilities, up to   inner
automorphisms of   vertex groups. But composing an inclusion with 
such an automorphism does not change the Bass-Serre tree [\Ba], so
we have shown that $\D$ only meets  finitely many $\Out(G)$-orbits of
non-metric reduced trees. Assertion (3) follows. 
\cqfd\enddemo

All finiteness conditions (but not smallness) are satisfied when $G$ is a
``generic'' hyperbolic group (one whose boundary is a Menger curve).  
When all conditions are satisfied,
$\Out(G)$ acts on the contractible complex
$P\D$ with finite stabilizers and finitely many orbits of simplices, so is 
$F_\infty$. For instance: 

\thm{Corollary \sta} Let $p$ be prime. Suppose  $T$ is an   
irreducible  locally finite  $G$-tree with stabilizers commensurable to
$BS(1,p)$. If there is no nontrivial integer modulus, then $\Out(G)$ is
$F_\infty$. 
\fthm

\demo{Proof} $BS(1,p)$ is small. Its subgroups of finite index are 
isomorphic to some
$BS(1,p^n)$, so have trivial center. Since $p$ is prime, their    Out is 
finite   by [\Co]. A group containing
$BS(1,p^n)$ with finite index also has Out finite  (see  e.g\. [\GL, Lemma
5.4]).
\cqfd \enddemo

\Refs\widestnumber\no{99}
\refno=0

\bref \by
 H. Bass \paper
Covering theory for graphs of groups \jour
J. Pure Appl. Alg. \vol89\pages3--47 \yr1993
\endref

\bref \by
 H. Bass, R.  Kulkarni\paper
Uniform tree lattices \jour
J. Am. Math. Soc. \vol3\pages843--902 \yr1990
\endref

\bref \by  M. Bestvina, M.  Feighn 
\paper Bounding the complexity of simplicial group actions on trees 
\jour Invent. Math. \vol103\pages 449--469 \yr1991 
\endref 

\bref \by B. Bowditch \paper Cut points and canonical splittings of
hyperbolic groups\jour Acta Math. \vol180\yr1998\pages145--186\endref

\bref\by M. Bridson, A. Miller \jour Lost manuscript
\endref

\bref \by I. Chiswell \book Introduction to
$\Lambda$-trees\publ World Scientific Publishing Co.\yr2001
\endref

\bref\by M. Clay \paper Contractibility of deformation spaces of
$G$-trees
\jour Alg. \& Geom. Topology\vol5\yr2005\pages1481--1503\endref

\bref\by M. Clay \paper A fixed point theorem for
deformation spaces of $G$-trees 
\jour Comm. Math. Helv. (to appear), arXiv:math.GR/0502248\endref

\bref\by M. Clay  \paper Deformation spaces of $G$-trees
\jour PHD thesis \yr2006\endref

\bref  \by D.J. Collins \paper The automorphism towers   of some
one-relator  groups\jour Arch. Math.\yr1983\vol40\pages385--400
\endref

\bref  \by M. Culler, J. Morgan \paper Group actions on \Rt s\jour
Proc. London Math. Soc.
\vol 55
\yr1987\pages571--604 \endref

\bref  \by M. Culler, K. Vogtmann \paper Moduli of graphs and
automorphisms of free groups\jour Invent. Math.
\yr1986\vol 84\pages 91--119
\endref

\bref  \by M. Culler, K. Vogtmann \paper The boundary of outer space in rank two,
{\rm in ``Arboreal group theory (Alperin, ed.)''} \jour MSRI Publ. 19\publ
Springer Verlag \yr1991\pages 189--230 
\endref

\bref\by W. Dicks\paper Groups, trees and projective modules
\jour Lecture Notes in Mathematics\vol 790
\yr1980
\endref

\bref \by M.J. Dunwoody, M.E. Sageev \paper JSJ-splittings for
finitely presented groups over slender groups\jour Invent. Math.
\vol135\yr1999\pages 25--44
\endref

\bref   \by M. Forester\paper Deformation and rigidity of simplicial
group actions on trees\jour Geom.
\& Topol.\vol 6\yr2002\pages 219--267\endref

\bref   \by M. Forester\paper On uniqueness of JSJ
decompositions of finitely generated groups\jour Comm. Math.
Helv. \vol78\yr2003\pages 740--751 \endref

\bref  \by M. Forester\paper Splittings of generalized
Baumslag-Solitar  groups \jour Geom\. Dedicata
(to appear), arXiv:math.GR/0502060
\endref

\bref \by K. Fujiwara, P. Papasoglu\paper  
JSJ-decompositions of
finitely presented groups and complexes of groups\jour
GAFA\vol16\yr2006\pages70--125
\endref

\bref \by V. Guirardel, G. Levitt \paper The outer space   of a free
product \jour Proc. London Math. Soc\. (to
appear), arXiv:math.GR/0501288\endref

\bref\by V. Guirardel, G. Levitt\paper A general construction of JSJ
decompositions\jour Proceedings of the 2005 Barcelona conference on
group theory, to appear, 
(available from http://www.picard.ups-tlse.fr/\~ {}guirardel/)  
\endref

\bref \by V. Guirardel, G. Levitt   \jour In preparation\endref

\bref \by S. Krsti\'c, K. Vogtmann\paper Equivariant outer space and
automorphisms of free-by-finite groups\jour Comm. Math. Helv. \yr
1993\vol 68\pages 216--262
\endref
\bref\by G. Levitt \paper Automorphisms of   hyperbolic
groups and graphs of groups \jour Geom.
Dedic. \vol 114\yr2005
\pages 49--70 
\endref

\bref\by G. Levitt\paper On the automorphism group of generalized
Baumslag-Solitar groups\hfill\break\jour arXiv:math.GR/0511083
\endref

      \bref   \by D. McCullough, A. Miller \paper Symmetric automorphisms
of free products \jour Mem. Amer. Math. Soc.
\vol122\yr1996\endref

\bref\by J.W. Morgan, P.B. Shalen \paper Valuations, trees, and
degenerations of hyperbolic structures, I \jour Ann.
Math\vol120\yr1984\pages401--476
\endref

\bref \by F. Paulin \paper The Gromov topology on $\R$-trees
\jour Topology and its App.\vol 32 \yr 1989\pages 197--221\endref

\bref \by F. Paulin \paper Sur les automorphismes ext\'erieurs des
groupes hyperboliques
\jour Ann. Scient. \'Ec. Norm. Sup. \vol 30  \yr1997\pages 147--167  
\endref

\bref  \by  E. Rips, Z.  Sela \paper
Cyclic splittings of finitely presented groups and the canonical  JSJ
decomposition\jour Ann. Math. \vol146\pages 55--109 \yr1997 
\endref

\bref
\by P. Scott, G.A. Swarup\paper Regular neighbourhoods and canonical
decompositions for groups\jour Ast\'erisque \vol289\yr2003
\endref

\bref \by Z. Sela\paper
Acylindrical accessibility for groups\jour
Invent. Math. \vol129\pages 527--565 \yr1997 
\endref 

\bref\by J.-P. Serre \paper Arbres, amalgames, SL$_2$\jour Ast\'erisque\vol 46
\yr1977
\endref

\bref   \by P.B. Shalen \paper Dendrology of groups: an introduction,
{\rm in ``Essays in group theory (S.M. Gersten, ed.)''} \jour MSRI Publ. 8\publ
Springer Verlag \yr1987 \endref

\bref \by R. Skora\paper  Deformations of length functions in
groups\jour
      Preprint
\yr1989\endref

\bref \by T. White \paper Fixed points of finite groups of free group
automorphisms\jour Proc. AMS \vol118\yr1993\pages 681--688
\endref

\endRefs

\address  V.G.: Laboratoire \'Emile Picard, umr cnrs 5580, Universit\'e
Paul Sabatier, 31062 Toulouse Cedex 4, France.\endaddress\email
guirardel\@math.ups-tlse.fr{}{}{}{}{}\endemail

\address  G.L.: LMNO, umr cnrs 6139, BP 5186, Universit\'e de Caen,
14032 Caen Cedex, France.\endaddress\email
levitt\@math.unicaen.fr{}{}{}{}{}\endemail

\enddocument